\documentclass[11pt,epsfig,letter]{article}
\usepackage{amssymb}
\usepackage{amsmath}
\usepackage{amsthm}
\usepackage{epsfig}
\usepackage[margin=1.15in]{geometry}

\newtheoremstyle{thmm}{1.5ex plus 1ex minus .2ex}{1.5ex plus 1ex minus
.2ex}{\rmfamily}{}{\bfseries}{}{1em}{}
\theoremstyle{thmm}
\newtheorem{theorem}{Theorem}[section]
\newtheorem{lemma}{Lemma}[section]
\newtheorem{proposition}{Proposition}[section]

\renewenvironment{proof}[1][Proof]{\noindent\textit{#1. } }{\hfill$\square$}
\allowdisplaybreaks 

\newcommand{\nn}{\nonumber}
\def \endproof{\vrule height8pt width 5pt depth 0pt}
\def\refe#1{(\ref{#1})}

\def\d{\delta}

\def\R{\mathbb{R}}

\def\U{{\bf U}}

\def\d{\,{\rm d}}

\def\u{{\bf u}}

\begin{document}
\title{\bf
Convergence of a decoupled mixed FEM for miscible
displacement in interfacial porous media}
\author{
Buyang Li \footnote{Department of Mathematics,
Nanjing University, Nanjing, 210093, China.
The work of the author was supported in part by NSF of China
(Grant No. 11301262) {\tt buyangli@nju.edu.cn} 
}
,~
Hongxing Rui \footnote{Department of Mathematics,
Shandong University, Jinan, China.  {\tt hxrui@sdu.edu.cn}
}
 ~~and~
Chaoxia Yang \footnote{Department of Mathematics,
China University of Petroleum, Qingdao, China.  {\tt yangcx@upc.edu.cn}
}
 }
\date{}

\maketitle

\begin{abstract}
In this paper, we study the stability and convergence 
of a decoupled and linearized mixed finite element method  (FEM)
for incompressible miscible displacement 
in a porous media whose permeability 
and porosity are discontinuous across some interfaces. 
We show that the proposed scheme has optimal-order 
convergence rate unconditionally, without restriction on the 
grid ratio (between the time-step size and spatial mesh size). 
Previous works all required 
certain restrictions on the grid ratio
except for the problem with globally
smooth permeability and porosity.
Our idea is to introduce an intermediate 
system of elliptic interface problems, 
whose solution is 
uniformly regular in each subdomain separated 
by the interfaces and its finite element solution 
coincides with the fully discrete solution of the 
original problem. In order to prove the boundedness
of the fully discrete solution, we study
the finite element discretization of
the intermediate system of elliptic interface problems. 
\end{abstract}

\section{Introduction}
\setcounter{equation}{0}

Numerical computation of miscible displacement in porous media has attracted much attention in recent decades due to its applications in reservoir
simulations and exploration of underground oil; see \cite{AE,ChenEwing,CZW, Dou,DFP,Peaceman}. The model describes the motion of a miscible fluid of two (or more) components in porous media, where the velocity of the fluid is given by Darcy's law
 \begin{align*}
 \u=-\frac{k(x)}{\mu(c)}\nabla p .
 \end{align*}
In the last equation, $p$ denotes the pressure of the fluid mixture, $k(x)$ denotes the permeability of the porous media, and $\mu(c)$ is the viscosity of the fluid depending on the concentration $c$ of the first component.
The incompressibility of the fluid is described by
\begin{align*}
 \nabla\cdot\u=q_I-q_P ,
 \end{align*}
 where $q_I$ and $q_P$ are given injection and production sources.
The concentration $c$ is governed by a convection-diffusion equation
\begin{align*}
\displaystyle\Phi(x)\frac{\partial c}{\partial t}-\nabla\cdot(D({\bf u},x)\nabla
c)+{\bf u}\cdot\nabla c= \hat c q_I-cq_I  .
\end{align*}
where $\Phi(x)$ denotes the porosity of the media and 
$D({\bf u},x)$ denotes the
diffusion-dispersion tensor, which is given by \cite{BB1,BB2}
\begin{align*}
D({\bf u},x) = \Phi(x)\bigg[d_0 I +F({\rm Pe})|\u|\left ( \alpha_1 I +
(\alpha_2 - \alpha_1) \frac{{\bf u} \otimes {\bf u}}{|{\bf u}|^2}
\right )\bigg] .
\end{align*}
In this formula,  $F({\rm Pe})= {\rm Pe}/({\rm Pe}+ d_{\rm r}) $ is a 
function of the local molecular Peclet number ${\rm Pe}=d_p |{\bf u}|$, 
where $d_0$, $\alpha_1$, $\alpha_2$, $d_r $ and $d_p$ are positive 
constants related to the porous media. It is straightforward to verify 
that 
$$
d_1|\xi|^2\leq D(\u,x)\xi\cdot\xi\leq (d_2+d_3|\u|)|\xi|^2,
\quad\forall~\xi\in\R^d ,
$$
for some positive constants $d_1$, $d_2$ and $d_3$.

Existence of weak or semiclassical solutions for the miscible 
displacement equations was studied in \cite{ChenEwing,Feng}, and
numerical analysis of the model has been done by many authors.
In particular, a Galerkin FEM was studied by Ewing and 
Wheeler \cite{EW}, and a Galerkin-mixed FEM was analyzed 
by Douglas et al \cite{DEW},  where the
Galerkin method was used to solve the parabolic 
concentration equation
and a mixed FEM was applied to solve the elliptic 
pressure equation. 
For both methods, a linearized
semi-implicit Euler scheme was used for the time stepping and 
optimal error estimates were presented roughly 
under the restriction $\tau =o(h)$. 
In \cite{Duran}, a characteristic method was applied to 
the parabolic concentration equation and the mixed 
FEM was used to solve the elliptic pressure equation. 
Optimal error estimates
were established under the same condition, i.e. $\tau=o(h)$.
More recently, a Galerkin method combined with a
post-process technique was studied in \cite{ML},
an Euler--Lagrange localized approximation 
method was studied in \cite{Wang} and 
a modified method of characteristics combined 
with mixed FEM was studied in \cite{SY}.
In all these works, error estimates were derived 
with certain restrictions on the grid ratio. 
To remove these restrictive conditions, a new approach 
was introduced in \cite{LS,LS2} to decouple the discretization 
errors from the temporal and spatial directions, and
optimal error estimates of a Galerkin-mixed FEM was 
established without restriction on the grid ratio. The 
methodology of \cite{LS,LS2} was later successfully applied 
to other nonlinear parabolic equations, such as the nonlinear 
Schr\"odinger equation \cite{JWang}, the thermistor equations 
\cite{HGao} and the Navier-Stokes equations \cite{SWS}.
However, all the analyses presented in these 
works rely on the global $H^2$ regularity of the 
``time-discrete solution'' (the solution of the linearized PDEs), 
which requires the permeability and porosity to be globally 
smooth in the miscible displacement model. 

In engineering computations, due to the existence of fault, 
filling-type karst caves or complex geological composition,  
the permeability and porosity are often discontinuous 
across some interfaces. It is desirable to solve the miscible 
displacement equations with discontinuous permeability and 
porosity by stable and accurate numerical methods.
For this purpose, numerical methods for flow in porous 
media with discontinuous permeability have been studied 
by many authors based on linear models. For example, 
see \cite{AH,CJMR,Ewing99,RWG} on the approximation 
of the elliptic pressure equation and see \cite{CZW} on the 
approximation of a parabolic pressure equation from the 
compressible model. Convergence of finite element methods 
for general linear elliptic and parabolic interface problems 
can also be found in \cite{CZ,LMWZ,SD05}. 
 
\begin{figure}[pth]
\centering
\includegraphics[width=1.8in]{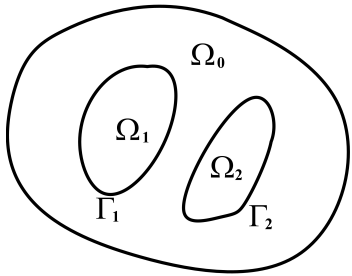}
\caption{The domain and 
the interfaces.}\label{Fig01}
\end{figure}

In this paper, we study stability and convergence of 
fully discrete FEMs for the full model of miscible 
displacement in porous media, where the permeability 
and porosity are discontinuous across some interfaces. 
Mathematically, we assume that the domain $\Omega$ 
is partitioned into $\Omega=\cup_{m=0}^{M}\Omega_m$ 
separated by the interfaces $\Gamma_m$, $m=1,\cdots,M$, 
as shown in Figure \ref{Fig01}, and we consider the nonlinear 
elliptic-parabolic interface problem
\begin{align}
&\left\{
\begin{array}{ll}
\displaystyle\Phi(x)\frac{\partial c}{\partial t}-\nabla\cdot(D({\bf u},x)\nabla
c)+{\bf u}\cdot\nabla c= \hat c q_I-cq_I &\mbox{in}~~\Omega_m,\\[10pt]
\left[c\right]=0,\quad \left[D(\u,x)\nabla c\cdot{\bf n}\right]=0 &\mbox{on}~~\Gamma_m,\\[10pt]
D(\u,x)\nabla c\cdot{\bf n}=0&\mbox{on}~~\partial\Omega,\\[8pt]
c(x,0)=c_0(x)~~ &\mbox{for}~~x\in\Omega ,
\end{array}
\right.\label{e-fuel-1}\\[8pt]
&\left\{
\begin{array}{ll}
\displaystyle  \nabla\cdot\u=q_I-q_P  &\mbox{in}~~\Omega_m,\\[8pt]
\displaystyle\u=-\frac{k(x)}{\mu(c)}\nabla p  &\mbox{in}~~\Omega_m,\\[10pt]
\displaystyle \left[ p \right]=0,\quad \left[\u\cdot{\bf n}\right]=0 &\mbox{on}~~\Gamma_m,\\[8pt]
\displaystyle \u\cdot{\bf n}=0    &\mbox{on}~~\partial\Omega .
\end{array}
\right.
\label{e-fuel-3}
\end{align}
In each subdomain $\Omega_m$, the pressure $p$,  the velocity $\u$
and the concentration $c$ are governed by the partial differential equations, and jump conditions are specified across the interfaces. The permeability $k(x)$ and porosity $\Phi(x)$ are assumed to be constant in each subdomain $\Omega_m$ but are discontinuous across the interfaces $\Gamma_m$. 

Clearly, the diffusion-dispersion tensor $D(\u,x)$ 
is an unbounded function of $\u$. Due to this strong 
nonlinearity and the coupling of equations, 
previous error estimates presented for the linear interface 
problems cannot be extended here. 
A direct application of the traditional error estimates requires 
undesired restrictions on the grid 
ratio to control the numerical velocity. In order to avoid any 
restrictive conditions on the grid ratio,
one has to use the error-splitting technique 
introduced in \cite{LS,LS2}.
However, due to the discontinuity of the 
permeability and porosity 
across the interfaces, the solution of 
\refe{e-fuel-1}-\refe{e-fuel-3} is not globally 
smooth. Instead, they are at most piecewise 
smooth \cite{BP70}, as assumed in this 
paper. In this case, the analysis for the Galerkin-mixed 
FEM presented in \cite{LS2} does not work.
In this paper, we show that a decoupled and linearized mixed FEM is 
stable for the nonlinear interface problem by proving 
that the time-discrete solution is piecewise 
smooth enough in each subdomain separated by the interfaces. 
Optimal error estimates are established 
without restriction on the grid ratio.
We believe that the methodology of this paper, together with 
Lemma \ref{LemHk0}-\ref{LemHkP1} introduced here, can also be applied to other
nonlinear parabolic interface problems in engineering and physics.


\section{Main results}
\setcounter{equation}{0}

Suppose that the smooth domain $\Omega$ is 
partitioned into  $\Omega=\cup_{m=0}^{M }\Omega_m$, 
where $\Omega_m$ is enclosed by
a smooth interface $\Gamma_m$ for $m=1,\cdots,M$, 
and $\Gamma_0=\partial\Omega$. 
For any integer $s\geq 0$ and a subdomain 
$\Omega_m$, we let $W^{s,p}(\Omega_m)$ and 
$H^s(\Omega_m):=W^{s,2}(\Omega_m)$ denote the
usual Sobolev spaces of functions defined on the 
domain $\Omega_m$; see \cite{Adams}. Let 
$L^p$ denote the abbreviations of $L^p(\Omega)$ 
and define $\overline W^{k,p}$ as the subspace of 
$L^p$ equipped with the norm
$$
\|f\|_{\overline W^{s,p}}:=\sum_{m=0}^{M }\|f\|_{W^{s,p}(\Omega_m)}  .
$$
Therefore, the functions in $\overline W^{s,p}$ are in 
$W^{s,p}(\Omega_m)$ for each subdomain $\Omega_m$, 
but may not be continuous in the whole domain $\Omega$.
To simplify the notations, 
we define $\overline H^s:=\overline W^{s,2}$,
 $\overline L^p:=\overline W^{0,p}$ and 
$$(f,g)=\sum_{m=0}^M\int_{\Omega_m}f(x)g(x)\d x,
\qquad\mbox{for}~~f,g\in \overline L^2 .$$ 
For any Banach space $X$ and a function 
$g:(0,T)\rightarrow X$, we define the norm
$$
\|g\|_{L^p((0,T);X)} =\left\{
\begin{array}{ll}
\displaystyle\biggl(\int_0^T\|g(t)\|_X^pdt\biggl)^\frac{1}{p}, &
1\leq p<\infty
,\\[10pt]
\displaystyle{\rm ess}\!\!\sup_{t\in (0,T)}\|g(t)\|_X, & p=\infty.
\end{array}
\right.
$$

Let $\{ t_n \}_{n=0}^N$ be a uniform partition of the time interval
$[0,T]$ with $\tau=T/N$ and denote
$$
p^n(x) = p(x,t_n),\quad {\bf u}^n(x)  = {\bf u}(x,t_n), \quad c^n(x)  = c(x,t_n) \, .
$$
For any sequence of functions $\{ f^n \}_{n=0}^N$, we define
$D_\tau f^{n+1}=(f^{n+1}-f^n)/\tau$.
Let $\pi_h$ denote a quasi-uniform partition 
of $\cup_{m=0}^M\Omega_m$ into triangles (or tetrahedrons)
$T_j$,
$j=1,\cdots,J$. For a triangle $T_j$ with two vertices 
on the boundary $\partial\Omega$ or an interface $\Gamma_m$,
we define $\widetilde T_j$ to be a 
triangle with one curved side which fit the 
boundary or the interfaces exactly, 
with the same vertices as $T_j$. 
Let  $h=\max_{1\leq j\leq
J}\{\mbox{diam}\,T_j\}$ denote the mesh size, and 
let ${\mathbb P}_r$ denote the space of 
polynomials of degree $r\geq 1$. 
We define the discontinuous finite
element space
\begin{align*}
&S_h^r =\{w_h\in L^2(\Omega): w_h|_{T_j}\in {\mathbb P}_r
\mbox{~for~each~element~$T_j\in \pi_h$ and~$\int_\Omega w_h
dx=0$}\} .
\end{align*}
Let $S_h^r(\Gamma_m)$ denote the space of 
functions in $S_h^r(\Omega_m)$ restricted to $\Gamma_m$.
To simplify the notations, we define $\Gamma_0=\partial\Omega$,
$\Gamma=\bigcup_{m=0}^M\Gamma_m$, 
and define $S_h^r(\Gamma)$ as the space of functions on 
$\Gamma$ whose restriction
to $\Gamma_m$ is in $S_h^r(\Gamma_m)$.
Let ${\bf H}_\Gamma^1$ be the space of vector-valued functions
${\bf v}\in (\overline H^1)^d$ such that 
$\nabla\cdot {\bf v}\in L^2$,
${\bf v}\cdot{\bf n}=0$ on $\partial\Omega$ and
$[{\bf v}\cdot{\bf n}]=0$ on $\Gamma_m$, $m=1,\cdots,M$. 
Let ${\bf H}_h^r$ denote the Raviart--Thomas mixed finite 
element subspace of ${\bf H}_\Gamma^1$ introduced in
\cite{DuranM,RT,VThomee}, which coincides with 
an element of ${\mathbb P}_r^d\oplus {\bf x}{\mathbb P}_r$ 
in each triangle $T_j$. Moreover, we require 
that the functions ${\bf v}_h\in {\bf H}_h^r$ satisfy the boundary 
condition
$\int_{\widetilde e_j}{\bf v}_h\cdot{\bf n}\, \chi_h\d s=0$, 
$\forall\chi_h\in S_h^r$, 
on each boundary edge $\widetilde e_j$ 
and the jump condition 
$\int_{\widetilde e_j}[{\bf v}_h\cdot{\bf n}]\, \chi_h\d s=0$, 
$\forall\chi_h\in S_h^r$, on each 
interface edge $\widetilde e_j$. 
Then we have $\nabla\cdot{\bf v}_h\in S_h^r$  
for ${\bf v}_h\in {\bf H}_h^r$.

To approximate $p$, $c$, $\u$ and ${\bf w}=-D(\u,x)\nabla c$, 
we look for $P_h^n,{\cal C}_h^n\in
S_h^r$ and $\U_h^n,{\bf W}_h^n\in {\bf H}_h^r$ which satisfy the equations
\begin{align}
& \Big(\frac{\mu( {\cal C}_h^{n})}{k(x)} \U_h^{n},\,{\bf v}_h \Big)
=\Big(P_h^{n} ,\, \nabla \cdot {\bf v}_h  \Big),
\label{e-FEM-1}\\[3pt]
& \Big(\nabla\cdot \U_h^{n} ,\, \varphi_h\Big) =\Big(q_I^{n}-q_P^{n},\,
\varphi_h\Big),
\label{e-FEM-2}\\[3pt]
&\Big(D( \U^{n}_h,x)^{-1}{\bf W}_h^{n+1},\,\overline {\bf v}_h\Big)
=\Big({\cal C}_h^{n+1} ,\, \nabla \cdot \overline {\bf v}_h  \Big),
\label{e-FEM-3}\\[3pt]
& \Big(\Phi(x) D_\tau {\cal C}_h^{n+1}, \,\overline\varphi_h\Big)
+ \Big(\nabla\cdot {\bf W}_h^{n+1}, \, \overline\varphi_h \Big)  
- \Big(D(\U^n_h,x)^{-1}\U^{n}_h\cdot {\bf W}_h^{n+1},\, 
\overline \varphi_h\Big) \nn\\
&\qquad\qquad\qquad\qquad\qquad\qquad\qquad\qquad = \Big(\hat
c^{n+1}q_I^{n+1}-{\cal C}_h^{n+1}q_I^{n+1}, \, \overline\varphi_h\Big),
 \label{e-FEM-4}
\end{align}
for any ${\bf v}_h,\overline {\bf v}_h\in {\bf H}_h^r$ 
and $\varphi_h,\overline\varphi_h\in S_h^r$,
where $n=0,1,2,\cdots$, and the initial 
data ${\cal C}_h^0$ is chosen as the
Lagrangian interpolation of $c^0$. 

For the initial-boundary value problem
(\ref{e-fuel-1})-(\ref{e-fuel-3}) to be 
well-posed, we require the compatibility condition
\begin{equation}\label{CompCond}
\int_\Omega q_I \d x=\int_\Omega q_P \d x ,
\end{equation}
and the physical restrictions
\begin{align}
&\| q_I \|_{L^\infty}+\| q_P \|_{L^\infty} \leq q_0  ,\label{qq}\\
&k_0^{-1}\leq k(x)\leq k_0\quad\mbox{for}~~x\in\Omega,\\
&\Phi_0^{-1}\leq \Phi(x)\leq \Phi_0\quad\mbox{for}~~x\in\Omega ,\\
&\mu_0^{-1}\leq \mu(s)\leq \mu_0\quad\mbox{and}\quad
|\mu'(s)|\leq \mu_0\quad\mbox{for}~~s\in\R ,
\end{align}
for some positive constants $q_0$, $k_0$, $\Phi_0$ and $\mu_0$.
Moreover, we assume that the solution of the
initial-boundary value problem (\ref{e-fuel-1})-(\ref{e-fuel-3}) exists
and possesses certain piecewise regularity such as
\begin{align}
\label{StrongSOlEST}
&\|p\|_{L^\infty((0,T);\overline H^{r+2} )}
+\|{\bf u}\|_{L^\infty((0,T);\overline H^{r+1} )}
+\|\partial_t{\bf u}\|_{L^\infty((0,T);\overline H^{1} )} 
+\|\partial_{tt}{\bf u}\|_{L^2((0,T);\overline L^{2} )} \nn \\
&+\|c\|_{L^\infty((0,T);\overline H^{r+2})}
+\|\partial_tc\|_{L^\infty((0,T);\overline H^{r+1} )} 
+\|\partial_{tt}c\|_{L^2((0,T);\overline L^2)}\leq C_0 ,
\end{align}
for some positive constant $C_0$.

The main result of this paper is the following theorem.
\medskip
\begin{theorem}\label{MainTHM}
{\it
Under the assumptions {\rm\refe{CompCond}-\refe{StrongSOlEST}}, 
there exists a positive constant $\tau_{**}$ such that when $\tau<\tau_{**}$ 
the finite element system
{\rm (\ref{e-FEM-1})-(\ref{e-FEM-4})} admits a unique solution $(P^n_h,
\U^n_h , {\cal C}^n_h,{\bf W}_h^n )$, $n=1,\cdots,N$, which satisfies that
\begin{align*}
&\max_{1\leq n\leq N}\big(\|P^n_h -p^n\|_{L^2}
+\| \U^n_h - \u^n\|_{L^2}
+\|{\cal C}^n_h - c^n\|_{L^2} \big) 
+\bigg(\sum_{n=1}^N\tau\|{\bf W}_h^n-{\bf w}^n\|_{L^2}^2\bigg)^\frac{1}{2}
\leq C_{**}(\tau+h^{r+1}),
\end{align*}
where $C_{**}$ is some positive constant 
independent of  $\tau$ and $h$.
}
\end{theorem}

The proof of Theorem \ref{MainTHM} is presented in 
Section \ref{SEction3}-\ref{SEction5}.
In Section 3, we introduce an intermediate problem, 
a system of elliptic interface problems,
whose finite element solution coincides with $(P^n_h,
\U^n_h , {\cal C}^n_h,{\bf W}_h^n )$, $n=0,1,2,\cdots$. 
Then we prove that the solution of the system of elliptic 
interface problems is piecewise smooth enough in each 
subdomain separated by the interfaces, and the piecewise 
regularity is uniform with respect to $\tau$ (as $\tau\rightarrow 0$). 
In Section 4, we present error estimates for the finite element 
discretization of the elliptic interface problems and prove the 
boundedness of the finite element solution based on the 
error estimates. In Section 5, we prove the error estimates in Theorem
\ref{MainTHM} based on the boundedness of the finite element solution.
Our analysis in Section 3
relies on the following two lemmas concerning 
the piecewise regularity of some elliptic and parabolic interface problems, 
which are generalizations of the results in \cite{BP70,CZ,DS,Huang} 
to problems with nonsmooth coefficients, 
with more precise dependence on the regularity of the coefficients. 
The proofs of the lemmas are given in Section \ref{LemHkP000}.

\begin{lemma}\label{LemHk0}
{\it If $A_{ij}\in \overline H^{2}$ satisfies that 
$K^{-1}|\xi|^2\leq \sum_{i,j=1}^dA_{ij}(x)\xi_i\xi_j\leq K|\xi|^2$ 
for $x\in\Omega$ and $\xi\in\R^d$, and $\phi\in H^1$ is a solution of
\begin{align}\label{HkEstTPEq}
\left\{
\begin{array}{ll}
\displaystyle
-\nabla\cdot\big(A\nabla \phi\big)=f &\mbox{in}~~\Omega,\\[10pt]
\displaystyle [\phi]=0,\quad \big[A\nabla \phi\cdot{\bf n}\big]=g_m 
&\mbox{on}~~\Gamma_m,~~m=1,\cdots,M,\\[10pt]
\displaystyle A\nabla \phi\cdot{\bf n}=0 &\mbox{on}~~\partial\Omega ,
\end{array}
\right. 
\end{align}
then
\begin{align}\label{HkEstTP}
\|\phi\|_{\overline H^k}\leq C\bigg(\|f\|_{\overline H^{k-2}} 
+\sum_{m=1}^M\|g_m\|_{H^{k-3/2}(\Gamma_m)}\bigg) \, ,\quad k=2,3 .
\end{align}
where the constant $C$ is independent of $\tau$.
 }
\end{lemma}

\begin{lemma}\label{LemHkP1}
{\it Suppose that 
\begin{align*}
&\displaystyle\max_{0\leq n\leq N-1}\sum_{i,j=1}^d
\big(\|A_{ij}^{n+1}\|_{\overline H^{2}}
+d_{n,0}\|D_\tau A^{n+1}_{ij}\|_{L^2}\big)\leq K,\\
& K^{-1}|\xi|^2\leq \sum_{i,j=1}^dA_{ij}^{n+1}(x)\xi_i\xi_j\leq K|\xi|^2,
\quad\forall~x\in\Omega,~ \xi\in\R^d ,
\end{align*}
and $\phi^{n+1}\in H^1$, 
$n=0,1,\cdots,N-1$, are solutions of
\begin{align}\label{HkEstPTEq}
\left\{
\begin{array}{ll}
\displaystyle
\Phi D_\tau\phi^{n+1}-\nabla\cdot\big(A^{n+1}\nabla \phi^{n+1}\big)=
f^{n+1}-\nabla\cdot{\bf g}^{n+1} &\mbox{in}~~\Omega,\\[10pt]
\displaystyle [\phi^{n+1}]=0,\quad \big[A^{n+1}\nabla \phi^{n+1}\cdot{\bf n}\big]=[{\bf g}^{n+1}\cdot{\bf n}] &\mbox{on}~~\Gamma_m,~~m=1,\cdots,M,\\[10pt]
\displaystyle A^{n+1}\nabla \phi^{n+1}\cdot{\bf n}=0 &\mbox{on}~~\partial\Omega ,\\[8pt]
\phi^0=0&\mbox{in}~~\Omega .
\end{array}
\right. 
\end{align}
Then we have
\begin{align}\label{HkEstPT2}
&\max_{0\leq n\leq m}\|\phi^{n+1}\|_{\overline H^1}^2
+\sum_{n=0}^m \tau\|\phi^{n+1}\|_{\overline H^2}^2  \\
&\leq C\|g^{m+1}\|_{L^2}^2
+C_\epsilon\sum_{n=0}^m\tau\big(\|f^{n+1}\|_{L^2}^2
 +\|{\bf g}^{n+1}\|_{\overline H^1}^2+\| \phi^{n+1}\|_{\overline H^1}^2 
) +\epsilon \sum_{n=0}^m\tau  \|D_\tau {\bf g}^{n+1}\|_{L^2}^2 d_{n,0} , 
\nn
\end{align}
where the constant $C_\epsilon$ 
$($dependent on $\epsilon$ $)$ is independent of $\tau$, and
$~d_{n,0}= 
\left\{
\begin{array}{ll}
0 ,&\mbox{if}~~n=0, \\
1 ,&\mbox{if}~~1\leq n\leq N-1 .
\end{array}
\right.
$
}
\end{lemma}

In the rest part of this paper, we denote by $C$ a generic positive constant and
by $\epsilon$ a small generic positive constant, which
are independent of $n$, $\tau$ and $h$.

\section{The linearized PDEs: a system of elliptic interface problems}\label{SEction3}
\setcounter{equation}{0}
We introduce $(P^n, {\cal C}^n)$, $n=0,1,2,\cdots$, as the solution of an iterative system of linear elliptic interface problems:
\begin{align}
&\left\{
\begin{array}{ll}
\displaystyle
-\nabla\cdot\bigg(\frac{k(x)}{\mu({\cal C}^{n})}\nabla P^{n}\bigg)
=q_I^{n}-q_P^{n}   &\mbox{in}~~\Omega_m,\\[8pt]
\displaystyle \left[P^{n}\right]=0,\quad \left[\frac{k(x)}{\mu({\cal C}^{n})}\nabla P^{n}\cdot{\bf n}\right]=0 &\mbox{on}~~\Gamma_m,\\[8pt]
\displaystyle \frac{k(x)}{\mu({\cal C}^{n})}\nabla P^{n}\cdot{\bf n}=0&\mbox{on}~~\partial\Omega,
\end{array}
\right.
\label{TDe-fuel-2}\\[10pt]
&\left\{
\begin{array}{ll}
\Phi(x) D_\tau {\cal C}^{n+1} -\nabla\cdot(D(\U^{n},x)\nabla {\cal
C}^{n+1}) + \U^{n} \cdot\nabla{\cal C}^{n+1}= \hat c^{n+1} q_I^{n+1}-{\cal
C}^{n+1}q_I^{n+1}   &\mbox{in}~~\Omega_m,\\[10pt]
\left[{\cal
C}^{n+1}\right]=0,\quad \left[D(\U^{n},x)\nabla {\cal
C}^{n+1}\cdot{\bf n}\right]=0 &\mbox{on}~~\Gamma_m,\\[10pt]
D(\U^{n},x)\nabla {\cal
C}^{n+1}\cdot{\bf n}=0&\mbox{on}~~\partial\Omega,
\end{array}
\right.
\label{TDe-fuel-3}
\end{align}
with the initial condition ${\cal C}^0=c_0$ and the 
normalization condition $\int_\Omega P^{n}dx =0 $. 
Existence and uniqueness 
of the solution for the 
linear elliptic interface problems 
\refe{TDe-fuel-2}-\refe{TDe-fuel-3} follow
iteratively, and it is easy to see that $P^0=p^0$ and ${\bf U}^0=\u^0$ 
at the initial time step. 
With this definition, the fully discrete solution $(P^n_h,
\U^n_h , {\cal C}^n_h,{\bf W}_h^n )$, $n=0,1,2,\cdots$, defined
in \refe{e-FEM-1}-\refe{e-FEM-4} 
can also be viewed as the 
finite element solution of 
\refe{TDe-fuel-2}-\refe{TDe-fuel-3}.

In this section, we establish the uniform piecewise regularity of 
$(P^n,{\cal C}^n)$ with respect to $\tau$. 
The following proposition is the main result of this section. 

\begin{proposition}\label{ErrestTDSol}
{\it
There exists a positive constant
$\tau_0$ such that when $\tau<\tau_0$, 
we have
\begin{align}
&\|P^{n}\|_{\overline H^3} +\|\U^n\|_{\overline H^2}+\|{\cal
C}^n\|_{\overline H^{3}}+\|D_\tau {\cal C}^n\|_{\overline H^1}\nn\\
&+\biggl(\sum_{n=1}^{N}\tau\|
D_\tau\U^{n}\|_{\overline H^2}^2\biggl)^{\frac{1}{2}}
+\biggl(\sum_{n=1}^{N}\tau\|
D_\tau {\cal C}^{n}\|_{\overline H^2}^2\biggl)^{\frac{1}{2}}
\leq C  . \label{StrongSOlESTTD}
\end{align}
}
\end{proposition}
The importance of this proposition is that the constant 
$C$ does not depend on $\tau$, which is the key 
to our error estimates in the next section.

\medskip
\noindent{\it Proof of Proposition {\rm\ref{ErrestTDSol}}}~~~ 
Let $e^n_p=P^n-p^n$, 
$e_c^n={\cal C}^n-c^n$ and 
$e_u^n={\bf U}^n-{\bf u}^n$. 
Comparing \refe{TDe-fuel-2}-\refe{TDe-fuel-3} 
with \refe{e-fuel-1}-\refe{e-fuel-3},
we see that  $e^n_p$, $e_c^n$ and $e_u^n$
satisfy the equations
\begin{align}\label{TDerr-fuel-1}
\left\{
\begin{array}{ll}
\displaystyle
-\nabla\cdot\biggl(\frac{k(x)}{\mu({\cal C}^{n})}\nabla e^{n}_p\biggl)
= \nabla\cdot\biggl[\biggl(\frac{k(x)}{\mu({\cal C}^{n})}
-\frac{k(x)}{\mu(c^{n})}\biggl)\nabla p^{n}\biggl],
&\mbox{in}~~\Omega_m,\\[10pt]
\displaystyle [e_p^{n}]=0,\quad -\bigg[\frac{k(x)}{\mu({\cal C}^{n})}\nabla e^{n}_p\cdot{\bf n}\bigg]=\bigg[\biggl(\frac{k(x)}{\mu({\cal C}^{n})}
-\frac{k(x)}{\mu(c^{n})}\biggl)\nabla p^{n}\cdot{\bf n}\bigg] &\mbox{on}~~\Gamma_m, \\[10pt]
\displaystyle \frac{k(x)}{\mu({\cal C}^{n})}\nabla e^{n}_p\cdot{\bf n}=0 &\mbox{on}~~\partial\Omega ,
\end{array}
\right. 
\end{align}
\begin{align}\label{TDerr-fuel-3}
\left\{
\begin{array}{ll}
\displaystyle
\Phi(x) D_\tau e_c^{n+1} -\nabla\cdot(D(\U^{n},x)\nabla e_c^{n+1}) \\
=-\U^{n}\cdot\nabla e_c^{n+1}
+\nabla\cdot\big((D(\U^{n},x)-D({\bf u}^{n+1},x))\nabla c^{n+1}\big)\\
~~~-(\U^{n}-{\bf u}^{n})\cdot\nabla c^{n+1}-e_c^{n+1}q_I^{n+1} +{\cal E}^{n+1},
&\mbox{in}~~\Omega_m,\\[10pt]
\displaystyle [e_c^{n+1}]=0,\quad \big[D(\U^{n},x)\nabla e_c^{n+1}\cdot{\bf n}\big]
=[(D(\u^{n+1},x)-D({\bf U}^{n},x))\nabla c^{n+1}\cdot{\bf n}] &\mbox{on}~~\Gamma_m, \\[10pt]
\displaystyle D(\U^{n},x)\nabla e_c^{n+1}\cdot{\bf n}=0 &\mbox{on}~~\partial\Omega ,\\[8pt]
e_c^0=0&\mbox{in}~~\Omega ,
\end{array}
\right. 
\end{align}
and
\begin{align}
& e_u^{n}=-\frac{k(x)}{\mu({\cal C}^{n})}\nabla e^{n}_p  -
\biggl(\frac{k(x)}{\mu({\cal C}^{n})}
-\frac{k(x)}{\mu(c^{n})}\biggl)\nabla p^{n},
\label{TDerr-fuel-2}
\end{align}
respectively, where
${\cal E}^{n+1} = \Phi (\partial_tc^{n+1}-D_\tau c^{n+1})
+({\bf u}^{n+1}-{\bf u}^n)\cdot\nabla c^{n+1}$ 
denotes the truncation error due to the time discretization. 
From the regularity assumption for $c$ in (\ref{StrongSOlEST}) we can
see that
\begin{align}\label{TDtruncerr}
&\|{\cal E}^{n+1}\|_{L^2}\leq C\tau^{1/2},\quad \sum_{n=0}^{N-1}\tau \|{\cal
E}^{n+1}\|_{L^2}^2 \leq C\tau^2  .
\end{align}

Integrating (\ref{TDerr-fuel-1}) against $e_p^{n}$, we get
\begin{align}
\|\nabla e_p^{n}\|_{L^2} &\leq \biggl
\|\biggl(\frac{k(x)}{\mu({\cal C}^{n} )}
-\frac{k(x)}{\mu(c^{n} )}\biggl)\nabla p^{n} \biggl \|_{L^2} \le C
\|e_c^{n} \|_{L^2}  \label{TDerr1fuel-1}
\end{align}
which together with (\ref{TDerr-fuel-2}) gives
\begin{align}
\|e_u^{n}\|_{L^2} & \leq C\|\nabla e_p^{n}\|_{L^2} +
C\|e_c^{n} \|_{L^2}  \leq  C
 \|e_c^{n} \|_{L^2} \, . \label{TDerr1fuel-2}
\end{align}
Then we integrate (\ref{TDerr-fuel-3}) against $e_c^{n+1}$ and obtain
\begin{align}
& D_\tau \biggl(\frac{1}{2} \|\sqrt{ \Phi}e_c^{n+1}\|_{L^2}^2 \biggl)
 + \| \sqrt{D(\U^{n},x)}\nabla e_c^{n+1} \|_{L^2}^2
\nn \\
& \leq C\| e^{n+1}_c\|_{L^2}^2 \|q_I^{n+1}-q_P^{n+1}\|_{L^\infty}
+C\big(\|e_u^{n}\|_{L^2} +\|u^{n+1}-u^n\|_{L^2} \big)\|
\nabla e_c^{n+1} \|_{L^2} \|\nabla c^{n+1}\|_{L^\infty}
\nn \\
&~~ + C\|e_u^{n}\|_{L^2} \|e_c^{n+1}\|_{L^2} \|\nabla
c^{n+1}\|_{L^\infty}+ C \|e_c^{n+1}\|_{L^2}^2 \| q_I^{n+1} \|_{L^\infty}  +
C\|{\cal E}^{n+1}\|_{L^2} \| e_c^{n+1} \|_{L^2}
\nn \\
& \leq \frac{1}{2} \|\nabla e_c^{n+1} \|_{L^2}^2 
+C\big(\|e_c^{n+1}\|_{L^2}^2 +\|e_c^{n}\|_{L^2}^2  
  + \|{\cal
E}^{n+1}\|_{L^2}^2+\|D_\tau\u\|_{L^2}^2\tau^2\big),
 \label{TDerr1-fuel-3} 
\end{align}
where we have used the inequality
$$
| ( \U^{n} \cdot \nabla e_c^{n+1}, \,  e_c^{n+1} )| =  | ( \nabla
\cdot \U^{n}, \, |e_c^{n+1}|^2 )| \leq  \| e^{n+1}_c \|_{L^2}^2
\|q_I^{n+1}-q_P^{n+1}\|_{L^\infty} .
$$
By applying Gronwall's inequality to \refe{TDerr1-fuel-3}, we see that 
there exists a positive constant $\tau_1$ such that when $\tau<\tau_1$ 
there holds
\begin{align}\label{fdsmkwjerqio984}
\max_{1\leq n\leq N}\|e_c^{n}\|_{L^2}^2 + \sum_{n=1}^{N}\tau \big \|\nabla e_c^{n}
\big \|_{L^2}^2 \leq C\tau^2 .
\end{align}
 The last
inequality, together with (\ref{TDerr1fuel-1})-(\ref{TDerr1fuel-2}), implies that
\begin{align}
&\max_{1\leq n\leq N}
\| D_\tau e_c^{n} \|_{L^2}^2 +
\sum_{n=1}^N\tau \| D_\tau e_c^{n} \|_{L^6}^2
\leq C, \label{fdnhifho08}\\
&\max_{1\leq n\leq N}
\big(\| e_u^{n} \|_{L^2} + \| \nabla e_p^{n} \|_{L^2} \big)
\leq C \tau .
\label{eupL2}
\end{align}

Let $d_{n,0}$ be the constant defined in Lemma \ref{LemHkP1}. 
We proceed with a mathematical induction on
\begin{align}\label{mathind1}
&
\|{\cal C}^{n}\|_{\overline H^2}+\|\U^{n}\|_{\overline H^2}
+d_{n,0}\|D_\tau \U^{n}\|_{L^2}\leq 
\max_{0\leq n\leq N}(\|c^{n}\|_{\overline H^2}
+\|\u^{n}\|_{\overline H^2}
+d_{n,0}\|D_\tau \u^{n}\|_{L^2})+1 ,
\end{align}
which clearly holds when $n=0$ 
(as ${\cal C}^0=c^0$ and $\U^0=\u^0$). 
We shall assume that the above inequality 
holds for $0\leq n\leq k$ and prove that it 
also holds for $n=k+1$.

With \refe{mathind1}, we can apply
Lemma \ref{LemHk0} to 
\refe{TDerr-fuel-1} for $0\leq n\leq k$ and
obtain
\begin{align*}
&\|e^{n}_p\|_{\overline H^2}
&\leq C\bigg\|\nabla\cdot\biggl[\biggl(\frac{k(x)}{\mu({\cal C}^{n})}
-\frac{k(x)}{\mu(c^{n})}\biggl)\nabla p^{n}\biggl]\bigg\|_{\overline L^2}
+C\bigg\|\left[\biggl(\frac{k(x)}{\mu({\cal C}^{n})}
-\frac{k(x)}{\mu(c^{n})}\biggl)\nabla p^{n}\right]\bigg\|_{H^{1/2}(\Gamma)}
\leq C\|e^{n}_c\|_{\overline H^1}  ,\\
&\|e^{n}_p\|_{\overline H^3}
&\leq C\bigg\|\nabla\cdot\biggl[\biggl(\frac{k(x)}{\mu({\cal C}^{n})}
-\frac{k(x)}{\mu(c^{n})}\biggl)\nabla p^{n}\biggl]\bigg\|_{\overline H^1}
+C\bigg\|\left[\biggl(\frac{k(x)}{\mu({\cal C}^{n})}
-\frac{k(x)}{\mu(c^{n})}\biggl)\nabla p^{n}\right]\bigg\|_{H^{3/2}(\Gamma)}
\leq C\|e^{n}_c\|_{\overline H^2}  ,
\end{align*}
and 
from \refe{TDerr-fuel-2} we see that
\begin{align}  \label{euW1p}
&\|e^{n}_u\|_{\overline H^1}
\leq C\|e^{n}_c\|_{\overline H^1} ,\\
&\|e^{n}_u\|_{\overline H^2}
\leq C\|e^{n}_c\|_{\overline H^2} . 
\end{align} 
As a consequence, by the Sobolev embedding inequality, we have
$$\|\nabla e^{n}_p\|_{L^\infty}\leq C\|e^{n}_p\|_{\overline H^3}\leq C
\qquad\mbox{and}\qquad 
\|e^{n}_u\|_{\overline H^2}\leq C ,
\qquad\mbox{for}~~0\leq n\leq k .$$ 
Applying the difference operator $D_\tau $ to the equation \refe{TDerr-fuel-3},
 we obtain
\begin{align}\label{EqDtauP}
\left\{
\begin{array}{ll}
\displaystyle
-\nabla\cdot\biggl(\frac{k(x)}{\mu({\cal C}^{n+1})}\nabla D_\tau e^{n+1}_p\biggl)
=\nabla\cdot{\bf f}^{n+1} ,
&\mbox{in}~~\Omega,\\[10pt]
\displaystyle [D_\tau e_p^{n+1}]=0,\quad 
-\bigg[\frac{k(x)}{\mu({\cal C}^{n+1})}\nabla 
D_\tau e^{n+1}_p\cdot{\bf n}\bigg]=\big[{\bf f}^{n+1}\cdot{\bf n}\big] &\mbox{on}~~\Gamma_m, \\[10pt]
\displaystyle \frac{k(x)}{\mu({\cal C}^{n+1})}
\nabla D_\tau  e^{n+1}_p\cdot{\bf n}=0 &\mbox{on}~~\partial\Omega ,
\end{array}
\right. 
\end{align}
where
\begin{align*}
{\bf f}^{n+1}=
D_\tau \bigg(\frac{k(x)}{\mu({\cal C}^{n+1})}\bigg)\nabla e^{n}_p
+D_\tau \biggl[\biggl(\frac{k(x)}{\mu({\cal C}^{n+1})}
-\frac{k(x)}{\mu(c^{n+1})}\biggl)\nabla p^{n+1}\biggl] .
\end{align*}
For $0\leq n\leq k$, from \refe{TDe-fuel-3} and \refe{TDerr-fuel-3} 
we can derive that
\begin{align*}
\|D_\tau {\cal C}^{n+1}\|_{L^2}&\leq C\|{\cal C}^{n+1}\|_{\overline H^2}
+C\|\U^{n}\|_{\overline H^2}\\
&\leq C\|e_c^{n+1}\|_{\overline H^2}+\|c^{n+1}\|_{\overline H^2}
+C\|\U^{n}\|_{\overline H^2}\\
& \leq 
C+C\|e_c^{n+1}\|_{\overline H^2} ,\\[5pt]
\|D_\tau e_c^{n+1}\|_{L^2}  &\leq C\|e_c^{n+1}\|_{\overline H^2}
+C\|e_u^{n}\|_{\overline H^2}+C\|{\cal E}^{n+1} \|_{L^2}\\
&\leq C\|e_c^{n+1}\|_{\overline H^2}
+C\|e_c^{n}\|_{\overline H^2}+C\|{\cal E}^{n+1} \|_{L^2},
\end{align*}
from \refe{EqDtauP} we see that
\begin{align*}
\|\nabla D_\tau e_p^{n+1}\|_{L^2}&\leq
C\|{\bf f}^{n+1}\|_{L^2}\\
&\leq C\|D_\tau {\cal C}^{n+1}\|_{L^2}\|\nabla e_p^n\|_{L^\infty}
+C\|D_\tau e_c^{n+1}\|_{L^2}+C\|e_c^{n+1}\|_{L^2}\|D_\tau c^{n}\|_{L^\infty}\\
&\leq C\|\nabla e_p^n\|_{L^\infty}+C\|e_c^{n+1}\|_{\overline H^2}\|\nabla e_p^n\|_{L^\infty}
+C\|e_c^{n+1}\|_{\overline H^2}
+C\|e_c^{n}\|_{\overline H^2}+C\|{\cal E}^{n+1} \|_{L^2}\\
&\leq C\|e_c^{n}\|_{\overline H^2}+C\|e_c^{n+1}\|_{\overline H^2}
+C\|{\cal E}^{n+1} \|_{L^2},
\end{align*}
and from \refe{TDerr-fuel-2} 
we derive that 
\begin{align} 
\|D_\tau e_u^{n+1}\|_{L^2} 
& \leq C(\|k(x)/\mu({\cal C}^{n+1})\|_{L^\infty}\|\nabla D_\tau e_p^{n+1}\|_{L^2} 
+\|D_\tau{\cal C}^{n+1}\|_{L^2}\|\nabla e_p^{n}\|_{L^\infty})\nn\\
&~~~+C(\|e_c^{n} \|_{L^\infty}\|D_\tau c^{n+1} \|_{L^2} 
+\|k(x)/\mu({\cal C}^{n+1})\|_{L^\infty}\|D_\tau e_c^{n+1} \|_{L^2} )
\|\nabla p^{n+1}\|_{L^\infty} 
 \nn\\
 &~~~+C\|e_c^{n+1} \|_{L^\infty}\|\nabla D_\tau p^{n+1}\|_{L^2}\nn\\
& \leq C(\|\nabla D_\tau e_p^{n+1}\|_{L^2} +
\|e_p^{n} \|_{\overline H^3} +\|e_c^{n+1}\|_{\overline H^2}
+\|e_c^{n}\|_{\overline H^2}
+\|D_\tau e_c^{n+1} \|_{L^2} ) \nn\\
&\leq C(
\|e_c^{n+1} \|_{\overline H^2} +\|e_c^{n} \|_{\overline H^2}
 +\|{\cal E}^{n+1} \|_{L^2} )  . \label{euW1p-2}
\end{align}
With \refe{mathind1}-\refe{euW1p-2}, we let
$$f^{n+1}=-\U^{n}\cdot\nabla e_c^{n+1}
+\nabla\cdot\big((D(\U^{n},x)-D({\bf u}^{n+1},x))\nabla c^{n+1}\big) 
-(\U^{n}-{\bf u}^{n})\cdot\nabla c^{n+1}-e_c^{n+1}q_I^{n+1} 
+{\cal E}^{n+1} ,
$$
and apply Lemma \ref{LemHkP1} to \refe{TDerr-fuel-3} 
for $0\leq n\leq k$. Then we derive that, for $0\leq m\leq k$,
\begin{align}\label{DtauU0}
&\max_{0\leq n\leq m}\|e^{n+1}_c\|_{\overline H^1}^2+\sum_{n=0}^m\tau\big\| e^{n+1}_c\big\|_{\overline H^2}^2 \nn  \\
&\leq C\max_{0\leq n\leq m}\|(D(\u^{n},x)-D({\bf u}^{n+1},x))\nabla c^{n+1}\|_{L^2}^2 \nn\\
&~~~ +\epsilon \sum_{n=0}^m \tau\|D_\tau((D(\U^{n},x)-D({\bf u}^{n+1},x))\nabla c^{n+1} )\|_{L^2}^2d_{n,0} \nn\\
&
~~~+C_\epsilon\sum_{n=0}^m\tau(\|f^{n+1}\|_{L^2}^2 +\|(D(\U^{n},x)-D({\bf u}^{n+1},x))\nabla c^{n+1}\|_{\overline H^1}+\|e^{n+1}_c\|_{\overline H^1}^2) \nn\\
&\leq C\tau^2+\epsilon \sum_{n=0}^m 
\tau(\|e_u^n\|_{L^\infty}^2+\|D_\tau e_u^n\|_{L^2}^2
+\|\u^{n+1}-\u^n\|_{L^\infty}^2+\|D_\tau (D(\u^{n+1},x)-D(\u^n,x))\|_{L^2}^2)d_{n,0}\nn\\
& ~~~+C_\epsilon\sum_{n=0}^m\bigg(
\|\U^{n}\|_{L^\infty}^2\|\nabla e_c^{n+1}\|_{L^2}^2
+\big(\|e_u^{n}\|_{ \overline H^1}^2+\tau^2\|D_\tau\u^{n+1}\|_{ \overline H^1}^2\big)
\|\nabla c^{n+1}\|_{L^\infty}^2 
 \nn\\
&~~~
+\big(\|e_u^{n}\|_{L^6}^2+\tau^2\|D_\tau\u^{n+1}\|_{ L^6}^2\big)
\|c^{n+1}\|_{\overline W^{2,3}}^2+\|e^{n}_u\|_{L^2}^2\|\nabla c^{n+1}\|_{L^\infty}^2
+\|e^{n+1}_c\|_{\overline H^1}^2+\|{\cal E}^{n+1}\|_{L^2}^2
\bigg) \nn\\
&\leq C\tau^2+\epsilon \sum_{n=0}^m 
\tau (\|e_c^n\|_{\overline H^2}^2 +\|e_c^{n-1}\|_{\overline H^2}^2 
+\|{\cal E}^{n}\|_{L^2}^2)d_{n,0} \\
&~~~
 +C_\epsilon\sum_{n=0}^m\tau \big(\|e^{n+1}_c\|_{\overline H^1}^2
+\|e_c^{n}\|_{\overline H^{1}}^2+\|{\cal E}^{n+1}\|_{L^2}^2
+\tau^2\|D_\tau\u^{n+1}\|_{ \overline H^1}^2\big) ,
\end{align}
which reduces to
\begin{align*}
&\max_{0\leq n\leq m}\|e^{n+1}_c\|_{\overline H^1}^2
+\sum_{n=0}^m\tau\big\| e^{n+1}_c
\big\|_{\overline H^2}^2 \leq C\tau^2 +\sum_{n=0}^m\tau\big( \|e^{n+1}_c\|_{\overline H^1}^2
+ \|e^{n}_c\|_{\overline H^1}^2) .
\end{align*}
By applying Gronwall's inequality, there exists a positive constant 
$\tau_3$ such that when $\tau<\tau_3$
\begin{align}\label{grdecH2}
\max_{1\leq n\leq k}\|e^{n+1}_c\|_{\overline H^1}^2
+\sum_{n=1}^k\tau\|e_c^{n+1}\|_{\overline H^{2}}^2  \leq C\tau^2.
\end{align}
From the last inequality we see that
\begin{align}\label{ecH222}
&\max_{0\leq n\leq k}\big(\|D_\tau {\cal C}^{n+1}\|_{\overline H^1}^2
+\|{\cal C}^{n+1}\|_{\overline H^2}^2\big)
+\sum_{n=0}^k\tau
\|D_\tau {\cal C}^{n+1}\|_{\overline H^{2}}^2
\leq C  ,\\
&\max_{0\leq n\leq k}\|e_c^{n+1}\|_{\overline H^{2}}
\leq C \tau^{1/2} , \nn
\end{align}
and from \refe{euW1p-2} we see that 
\begin{align*} 
\|D_\tau e_u^{k+1}\|_{L^2} 
\leq C(
\|e_c^{k+1} \|_{\overline H^2} +\|e_c^{k} \|_{\overline H^2}  +\|{\cal E}^{k+1} \|_{L^2} )
\leq C\tau^{1/2}  . 
\end{align*}
With $\displaystyle\max_{0\leq n\leq k}\|{\cal C}^{n+1}\|_{\overline H^2}\leq C$ 
given by \refe{ecH222}, 
we can apply
Lemma \ref{LemHk0} to \refe{TDerr-fuel-1} again and obtain
\begin{align*} 
&\|e^{k+1}_p\|_{\overline H^3}+\|e^{k+1}_u\|_{\overline H^2}
\leq C\|e^{k+1}_c\|_{\overline H^2}\leq C\tau^{1/2}  , \\
&\sum_{n=0}^k\tau \|e_u^{n+1}\|_{\overline H^2}^2
\leq C\sum_{n=0}^k\tau \|e_c^{n+1}\|_{\overline H^2}^2
\leq C\tau^2 ,
\end{align*} 
and so
\begin{align} 
&\max_{0\leq n\leq k}(\|P^{n+1}\|_{\overline H^3}^2
+\|\U^{n+1}\|_{\overline H^2}^2)
+\sum_{n=0}^k\tau \|D_\tau\U^{n+1}\|_{\overline H^2}^2\leq C,
\label{ecuH222}
\end{align} 
The last five inequalities imply that
there exists a positive constant $\tau_4$ such that when $\tau<\tau_4$ we have
\begin{align*} 
&
\|e_c^{k+1}\|_{\overline H^2}+\|e_u^{k+1}\|_{\overline H^2}
+\|D_\tau e_u^{k+1}\|_{L^2}\leq 
 1 .
\end{align*}
The mathematical induction on \refe{mathind1} is completed. 
Thus \refe{ecH222}-\refe{ecuH222} hold for $k=N-1$ with the 
same constant $C$, provided $\tau<\tau_0:=\min(\tau_1,\tau_2,\tau_3,\tau_4)$.

With the regularity $\displaystyle\max_{0\leq n\leq N}
\|D_\tau{\cal C}^n\|_{\overline H^1}\leq C$, as shown in \refe{ecH222}, 
by applying Lemma \ref{LemHk0} to \refe{TDe-fuel-3} we obtain
\begin{align}\label{CH3Est}
&\max_{1\leq n\leq N}\|{\cal C}^{n}\|_{\overline H^3}\leq C  .
\end{align}

The proof of Proposition \ref{ErrestTDSol} is completed. ~\endproof

\section{Boundedness of the $\U^n_h$}
\label{SEction4}
\setcounter{equation}{0}

Based on the finite element discretization of the elliptic interface problems, we prove the following proposition in this section.
\begin{proposition}\label{BdUnh}
{\it
There exist positive constants $\tau_*$ and $h_*$ such that when $\tau<\tau_*$ and $h<h_*$, the finite element system
{\rm (\ref{e-FEM-1})-(\ref{e-FEM-4})} admits a unique solution $(P^n_h,
\U^n_h,{\bf W}^n_h , {\cal C}^n_h )$, $n=1,\cdots,N$, such that
\begin{align*}
&\max_{1\leq n\leq N} \|\U^n_h \|_{L^\infty}  \leq C .
\end{align*}
}
\end{proposition}

Before we prove this proposition, we define some notations below.
Let $L_h$ denote the piecewise linear Lagrange interpolation operator onto the finite element space $S_h^r$.
Let $\Pi_h$ denote the $L^2$ projection onto the finite element space $S_h^r$, i.e.
$$
(\phi-\Pi_h \phi,\chi_h)=0,~~~\forall~\phi\in
L^2,~\chi_h\in S_h^r ，
$$
and let $\Pi^\Gamma_{h}$ denote the $L^2$ projection onto the finite element space $S_h^r(\Gamma)$ satisfying
$$
(\phi-\Pi^\Gamma_{h} \phi,\chi_h)_{\Gamma}=0,~~~\forall~\phi\in
L^2(\Gamma),~\chi_h\in S_h^r(\Gamma).
$$

Let $Q_h: {\bf H}_\Gamma^1 \rightarrow {\bf H}_h^r$ be a
projection satisfying (see \cite{DuranM,VThomee} for the construction of such a projection operator)
\begin{align}
&\big(\nabla\cdot({\bf v}-Q_h{\bf v})\, , \chi_h \big) = 0 ,\qquad~\,
\forall~\chi_h\in S_h^r, ~ {\bf v}\in  {\bf H}_\Gamma^1 , \label{Q}\\
&\int_{\widetilde e_j}({\bf v}-Q_h{\bf v})\cdot{\bf n}\, \chi_h \d s= 0 ,\quad ~
\forall~\chi_h\in S_h^r, ~ {\bf v}\in {\bf H}_\Gamma^1 ,
\end{align}
for any edge $\widetilde e_j$ in the triangulation. Then we have
$$
\|{\bf v}-Q_h {\bf v}\|_{L^2} 
+\|{\bf v}-Q_h {\bf v}\|_{L^2(\Gamma_m)} h^{1/2}
+\|\nabla\cdot({\bf v}-Q_h {\bf v})\|_{L^2} h
\leq C\|{\bf v}\|_{\overline H^k}h^k,\quad k=1,2,\cdots 
$$
Let ${\bf W}^{n+1}=-D(\U^{n},x)\nabla {\cal C}^{n+1}\in {\bf H}^1_\Gamma$ and, 
for any fixed integer $n\geq -1$, let 
$(\overline {\cal C}^{n+1}_h,\overline {\bf W}^{n+1}_h)
\in
S_h^r\times 
{\bf H}_h^r$ be the finite element solution of the equation
\begin{align}\label{MEPrjc}
\left\{
\begin{array}{ll}
(\nabla\cdot (\overline {\bf W}_h^{n+1}-{\bf W}^{n+1}),\chi_h)=0, &\forall~\chi_h\in S_h^r, \\[5pt]
(D(\U^{n},x)^{-1}(\overline {\bf W}_h^{n+1}-{\bf W}^{n+1}),{\bf v}_h) =(\overline{\cal C}^{n+1}_h-{\cal C}^{n+1},\nabla\cdot {\bf v}_h),&\forall~ {\bf v}_h\in {\bf H}_h^r ,
\end{array}
\right.
\end{align}
with $\int_\Omega(\overline{\cal C}^{n+1}_h-{\cal C}^{n+1})d x=0$ for the uniqueness of solution, where we define $\U^{-1}:=\U^{0}$. The pair $(\overline {\cal C}^{n+1}_h,\overline {\bf W}^{n+1}_h)$ can be viewed as the Ritz projection of $({\cal C}^{n+1},{\bf W}^{n+1})$ by the mixed FEM.

We require $\tau<\tau_0$ so that Proposition \ref{ErrestTDSol} holds.  With the regularity of ${\cal C}^{n+1}$ and ${\bf W}^{n+1}$ given in Proposition \ref{ErrestTDSol}, by the theory of mixed FEM for
linear elliptic equations \cite{DuranM,VThomee}, we have
\begin{align} 
&\|\phi-\Pi_h\phi\|_{L^2}\leq
Ch^2\|\phi\|_{\overline H^2}, \quad \forall~\phi\in \overline H^2,  
 \label{p-app0}\\[5pt]
 &\|\phi-\Pi^\Gamma_{h}\phi\|_{L^2(\Gamma)}\leq
Ch^{k+1/2}\|\phi\|_{H^{k+1/2}(\Gamma)}, \quad \forall~\phi\in \overline H^{k+1/2}(\Gamma)~~\mbox{with $1\leq k\leq r$},
 \label{p-apap}\\[5pt]
&\| {\bf W}^{n+1}-L_h{\bf W}^{n+1}\|_{L^2}
+\|{\bf W}^{n+1}-\overline {\bf W}_h^{n+1}\|_{L^2}  \leq
Ch^2, \\[5pt]
&\|{\cal C}^{n+1}-\overline {\cal C}_h^{n+1}\|_{L^2}\leq
Ch^2 . \label{p-app}
\end{align}
Therefore, by the inverse inequality, we have
\begin{align*}
&\|Q_h {\bf W}^{n+1}-L_h{\bf W}^{n+1}\|_{L^\infty}
+\|\overline{\bf W}^{n+1}-L_h{\bf W}^{n+1}\|_{L^\infty}\\
&\leq Ch^{-d/2}\|Q_h {\bf W}^{n+1} -L_h{\bf W}^{n+1}\|_{L^2}
\leq Ch^{2-d/2} ,
\end{align*}
which implies the existence of a positive constant $h_1$ such that when $h<h_1$ there holds
\begin{align} \label{p-app02}
\|Q_h{\bf W}^{n+1}\|_{L^\infty}+\|\overline{\bf W}^{n+1}\|_{L^\infty}
& \leq 2\|L_h{\bf W}^{n+1}\|_{L^\infty} +1 
\leq C .
\end{align}
Moveover, we need the following two lemmas in the proof of Proposition \ref{BdUnh}.

\begin{lemma}\label{DtauCh}
{\it Under the regularity of ${\cal C}^{n+1}$ and $\U^{n+1}$ proved in Proposition {\rm\ref{ErrestTDSol}}, we have\begin{align}
&\bigg(\sum_{n=0}^{N-1}\tau\|D_\tau ({\cal C}^{n+1}-\overline {\cal C}_h^{n+1})\|_{L^2}^2\bigg)^{1/2}\leq Ch^2 .   
\end{align}
}
\end{lemma}
\noindent{\it Proof}~~~
From \refe{MEPrjc} we derive that
\begin{align}\label{MEPrjc22aa}
&(\nabla\cdot(D_\tau\overline {\bf W}^{n+1}_h-D_\tau{\bf W}^{n+1})),\chi_h)=0 ,
\quad
\forall~ \chi_h\in S_h^r ,\\
&(D(\U^{n},x)^{-1}(D_\tau\overline {\bf W}^{n+1}_h-D_\tau{\bf W}^{n+1}),{\bf v}_h)+(D_\tau D(\U^{n},x)^{-1}(\overline {\bf W}^{n}_h-{\bf W}^{n}),{\bf v}_h) \nn\\
&=(D_\tau(\overline{\cal C}^{n+1}_h-\Pi_h{\cal C}^{n+1}),\nabla\cdot {\bf v}_h),\quad
\forall~ {\bf v}_h\in {\bf H}_h^r ,
\label{MEPrjc22bb}
\end{align}
where \refe{MEPrjc22aa} implies that 
$\nabla\cdot(D_\tau\overline {\bf W}^{n+1}_h-Q_hD_\tau{\bf W}^{n+1})=0$.
By choosing ${\bf v}_h=D_\tau\overline {\bf W}_h^{n+1}-Q_hD_\tau{\bf W}^{n+1}$ in 
\refe{MEPrjc22bb}, we derive that
\begin{align*} 
\|D_\tau\overline {\bf W}_h^{n+1}-Q_hD_\tau{\bf W}^{n+1}\|_{L^2}
&\leq C\|D_\tau {\bf W}^{n+1}-Q_hD_\tau {\bf W}^{n+1}\|_{L^2}
+C\|{\bf W}^{n+1}-\overline{\bf W}_h^{n+1}\|_{L^2} ,
\end{align*}
and so, by the inverse inequality, 
\begin{align*} 
&\|\nabla\cdot(D_\tau\overline {\bf W}_h^{n+1}-Q_hD_\tau{\bf W}^{n+1})\|_{L^2}
\leq Ch^{-1}\|D_\tau\overline {\bf W}_h^{n+1}-Q_hD_\tau{\bf W}^{n+1}\|_{L^2}
\leq C\|D_\tau {\bf W}^{n+1}\|_{H^1}+C .
\end{align*}

Let ${\bf v}=-D(\U^{n},x)\nabla g$, where $g$ is the solution of the elliptic interface problem
\begin{align*}
\left\{\begin{array}{ll}
-\nabla\cdot\big(D(\U^{n},x)\nabla g\big)=D_\tau(\overline{\cal C}^{n+1}_h-\Pi_h{\cal C}^{n+1})
&\mbox{in}~~\Omega_m,\\
\left[g\right]=0,\quad
[D(\U^{n},x)\nabla g\cdot{\bf n}]=0 &
\mbox{on}~~\Gamma_m,\\
D(\U^{n},x)\nabla g\cdot{\bf n}=0 &\mbox{on}~~\partial\Omega .
\end{array}
\right.
\end{align*}
Substituting ${\bf v}_h=Q_h{\bf v}$ into \refe{MEPrjc22bb}, we obtain
\begin{align*} 
&\|D_\tau ({\cal C}^{n+1}-\overline {\cal C}_h^{n+1})\|_{L^2}^2\\
&=(D(\U^{n},x)^{-1}D_\tau(\overline {\bf W}^{n+1}_h-{\bf W}^{n+1}),{\bf v})\\
&~~~ +(D(\U^{n},x)^{-1}D_\tau(\overline {\bf W}^{n+1}_h-{\bf W}^{n+1}),{\bf v}_h-{\bf v})
+(D_\tau D(\U^{n},x)^{-1}(\overline {\bf W}^{n}_h-{\bf W}^{n}),{\bf v}_h)\\
&=(\nabla\cdot (D_\tau \overline{\bf W}_h^{n+1}-D_\tau {\bf W}^{n+1}),g)
-([(D_\tau \overline{\bf W}_h^{n+1}-D_\tau {\bf W}^{n+1})\cdot{\bf n}],g)_{\Gamma}\\
&~~~ +(D(\U^{n},x)^{-1}D_\tau(\overline {\bf W}_h^{n+1}-{\bf W}^{n+1}),{\bf v}_h-{\bf v})
+(D_\tau D(\U^{n},x)^{-1}(\overline  {\bf W}_h^{n}-{\bf W}^{n}),{\bf v}_h)\\
&=(\nabla\cdot (D_\tau \overline{\bf W}_h^{n+1}-D_\tau {\bf W}^{n+1}),g-L_hg)
-([(D_\tau \overline{\bf W}_h^{n+1}-D_\tau {\bf W}^{n+1})\cdot{\bf n}],g-\Pi^\Gamma_{h}g)_{\Gamma}\\
&~~~ +(D(\U^{n},x)^{-1}D_\tau(\overline {\bf W}_h^{n+1}-{\bf W}^{n+1}),{\bf v}_h-{\bf v})
+(D_\tau D(\U^{n},x)^{-1}(\overline  {\bf W}_h^{n}-{\bf W}^{n}),{\bf v}_h)\\
&\leq C\|D_\tau{\bf W}^{n+1}\|_{\overline H^1}\|g-L_hg\|_{L^2}
+C\|[(D_\tau\overline  {\bf W}_h^{n+1}-D_\tau{\bf W}^{n+1})\cdot{\bf n}]
\|_{L^2(\Gamma)}\|g-\Pi^\Gamma_{h}g\|_{L^2(\Gamma)}\\
&~~~ +\|D_\tau\overline {\bf W}_h^{n+1}-D_\tau{\bf W}^{n+1} \|_{L^2}\|{\bf v}_h-{\bf v}\|_{L^2}
+\|D_\tau D(\U^{n},x)^{-1}\|_{L^\infty}
\|\overline {\bf W}^{n}_h-{\bf W}^{n}\|_{L^2}\|{\bf v}_h\|_{L^2}\\
&\leq C\|D_\tau{\bf W}^{n+1}\|_{\overline H^1}\|g\|_{\overline H^2}h^2
+C\|[(D_\tau\overline {\bf W}_h^{n+1}-Q_hD_\tau{\bf W}^{n+1})\cdot{\bf n}]\|_{L^2(\Gamma)}\|g\|_{H^{3/2}(\Gamma)}h^{3/2}\\
&~~~ +C\|[(Q_hD_\tau\overline {\bf W}^{n+1}-D_\tau{\bf W}^{n+1})\cdot{\bf n}]\|_{L^2(\Gamma)}\|g\|_{H^{3/2}(\Gamma)}h^{3/2}\\
&~~~ +C\|D_\tau {\bf W}^{n+1}\|_{\overline H^1}\|g\|_{\overline H^2}h^2
+C \|D_\tau \U^{n}\|_{L^\infty}\|{\bf v}_h\|_{L^2}\\
&\leq C\|D_\tau{\bf W}^{n+1}\|_{\overline H^1}\|g\|_{\overline H^2}h^2
+Ch^{-1/2}\|D_\tau\overline {\bf W}_h^{n+1}-Q_hD_\tau{\bf W}^{n+1}
\|_{L^2}\|g\|_{H^{3/2}(\Gamma)}h^{3/2}\\
&~~~
+C\|D_\tau{\bf W}^{n+1}\|_{\overline H^1}\|g\|_{H^{3/2}(\Gamma)}h^{2}+C\|D_\tau \U^{n}\|_{\overline H^2}\|g\|_{\overline H^2}h^2\\
&\leq C\|D_\tau{\bf W}^{n+1}\|_{\overline H^1}\|g\|_{\overline H^2}h^2
+C\|D_\tau \U^{n}\|_{\overline H^2}\|g\|_{\overline H^2}h^2\\
&\leq C(\|D_\tau{\cal C}^{n+1}\|_{\overline H^2}+\|D_\tau \U^{n}\|_{\overline H^2}\|)\|D_\tau ({\cal C}^{n+1}-\overline {\cal C}_h^{n+1})\|_{L^2} h^2,
\end{align*}
which reduces to
\begin{align*} 
&\|D_\tau ({\cal C}^{n+1}-\overline {\cal C}_h^{n+1})\|_{L^2}
\leq C(\|D_\tau{\cal C}^{n+1}\|_{\overline H^2}+\|D_\tau \U^{n}\|_{\overline H^2})h^2 .
\end{align*}
The last inequality, together with Proposition \ref{ErrestTDSol}, 
gives \refe{DtauCh}. ~\endproof\bigskip

\begin{lemma}\label{Lembf41}
{\it If $g\in \overline H^{k+1}$ with $1\leq k\leq r$, 
then
\begin{align*} 
\big|\big(g-\Pi^\Gamma_{h}g,{\bf v}_h\cdot{\bf n}\big)_{\partial\Omega}\big|+\sum_{m=1}^M\big|\big(g-\Pi^\Gamma_{h}g,[{\bf v}_h\cdot{\bf n}]\big)_{\Gamma_m}\big|\leq C\|g\|_{\overline H^{k+1}}\|{\bf v}_h\|_{L^2}h^{k+1} .
\end{align*}
}
\end{lemma} 
\noindent{\it Proof}~~~ For simplicity, we only prove the 2D case. The 3D case can be proved in the same way. For a triangle $T_j$ on the boundary (or an interface), we denote by $e_j$ its edge with at two vertices on the boundary (or an interface) and denote by $\widetilde e_j$ the curved edge on the boundary (or an interface). Via a rigid rotation, we assume that $e_j$ is on the $x_1$-axis and $\widetilde e_j$ is parametrized by $(x_1,y(x_1))$. Let $\chi_h\in S_h^r$ be a finite element function whose restriction to $e_j$ coincides with ${\bf v}_h\cdot{\bf n}$, i.e.
${\bf v}_h(x_1,0)\cdot{\bf n}(x_1,0)=\chi_h(x_1,0)$, satisfying $\|\chi_h\|_{L^2(T_j)}\leq C\|{\bf v}_h\|_{L^2(T_j)}$. Then we have, with $\d s=\sqrt{1+|y'(x_1)|^2}\d x_1$,
\begin{align*} 
&\int_{\widetilde e_j}(g-\Pi^\Gamma_{h}g){\bf v}_h\cdot{\bf n}\d s \\
&=\int_{e_j}\big(g(x_1,y(x_1))
-(\Pi^\Gamma_{h}g)(x_1,y(x_1))\big)
{\bf v}_h(x_1,y(x_1))\cdot{\bf n}(x_1,y(x_1)) \d s  \\
&=\int_{e_j}\big(g(x_1,y(x_1))-(\Pi^\Gamma_{h}g)(x_1,y(x_1))\big)
\big({\bf v}_h(x_1,y(x_1))\cdot{\bf n}(x_1,y(x_1)) 
-{\bf v}_h(x_1,0)\cdot{\bf n}(x_1,0) \big) \d s\\
&~~~ +\int_{e_j}\big(g(x_1,y(x_1))-(\Pi^\Gamma_{h}g)(x_1,y(x_1))\big)
(\chi_h(x_1,0)-\chi_h(x_1,y(x_1)) \d s\\
&\leq C\|g-\Pi^\Gamma_{h}g\|_{L^2(e_j)}\big(
\|\sup_{x_2\in T_j}\partial_y{\bf v}_h(\cdot,x_2)\|_{L^2(e_j)}h^2
+\|{\bf v}_h(\cdot,0)\|_{L^2(e_j)}h
+\|\sup_{x_2\in T_j}\partial_y\chi_h(\cdot,x_2)\|_{L^2(e_j)}h^2\big)\\
&\leq C\|g\|_{H^{k+1/2}(e_j)}(\|{\bf v}_h\|_{L^2(T_j)}
+\|\chi_h(\cdot,x_2)\|_{L^2(T_j)})h^{k+1}\\
&\leq C\|g\|_{H^{k+1}(T_j)}\|{\bf v}_h\|_{L^2(T_j)}h^{k+1} ,
\end{align*}
which implies that
\begin{align*} 
\int_{\partial\Omega}(g-\Pi^\Gamma_{h}g){\bf v}_h\cdot{\bf n}\d s
&\leq C\sum_{j}\|g\|_{H^{k+1}(T_j)}\|{\bf v}_h\|_{L^2(T_j)}h^{k+1}
\leq C\|g\|_{\overline H^{k+1}}\|{\bf v}_h\|_{L^2}h^{k+1}.
\end{align*}

The estimate of $\int_{\Gamma_m}(g-\Pi^\Gamma_{h}g)
[{\bf v}_h\cdot{\bf n}]\d s $ on an interface $\Gamma_m$
is similar.
~\endproof\bigskip\medskip

\noindent{\it Proof of Proposition {\rm\ref{BdUnh}}}~~~ Let $\tau<\tau_0$ so that \refe{p-app0}-\refe{p-app02} hold. The mixed weak formulation of \refe{TDe-fuel-2}-\refe{TDe-fuel-3} is
\begin{align}
& \Big(\frac{\mu( {\cal C}^{n})}{k(x)} \U^{n},\,{\bf v}_h \Big)
=\Big(P^{n} ,\, \nabla \cdot {\bf v}_h  \Big) -\Big(P^n,{\bf v}_h\cdot{\bf n}\Big)_{\partial\Omega}-\sum_{m=1}^M\Big(P^n,[{\bf v}_h\cdot{\bf n}]\Big)_{\Gamma_m},
\label{tde-FEM-1}\\[3pt]
& \Big(\nabla\cdot \U^{n} ,\, \varphi_h\Big) =\Big(q_I^{n}-q_P^{n},\,
\varphi_h\Big),
\label{tde-FEM-2}\\[3pt]
& \Big(D( \U^{n},x)^{-1}{\bf W}^{n+1},\,\overline {\bf v}_h\Big)
=\Big({\cal C}^{n+1} ,\, \nabla \cdot \overline {\bf v}_h  \Big)
-\Big({\cal C}^{n+1},\overline {\bf v}_h\cdot{\bf n}\Big)_{\partial\Omega}
-\sum_{m=1}^M\Big({\cal C}^{n+1},[\overline {\bf v}_h\cdot{\bf n}]\Big)_{\Gamma_m},
\label{tde-FEM-3}\\[3pt]
& \Big(\Phi(x) D_\tau {\cal C}^{n+1}, \,\overline\varphi_h\Big)
+ \Big(\nabla\cdot {\bf W}^{n+1}, \, \overline\varphi_h \Big)  - \Big(D(\U^n,x)^{-1}\U^{n}\cdot {\bf W}^{n+1},\, \overline \varphi_h\Big) \nn\\
&\qquad\qquad\qquad\qquad\qquad\qquad\qquad\qquad = \Big(\hat
{\cal C}^{n+1}q_I^{n+1}-{\cal C}^{n+1} q_I^{n+1}, \, \overline\varphi_h\Big),
 \label{tde-FEM-4}
\end{align}
for any ${\bf v}_h,\overline {\bf v}_h\in {\bf H}_h^r$ and $\varphi_h,\overline\varphi_h\in S_h$. The
above equations with the finite element system
(\ref{e-FEM-1})-(\ref{e-FEM-4}) imply that
\begin{align}
& \biggl(\frac{\mu({\cal C}^{n}_h)}{k(x)} \U^{n}_h-\frac{\mu({\cal C}^{n})}{k(x)} \U^{n},\,{\bf v}_h\biggl) =\Big(P^{n}_h-\Pi_hP^{n},\,
\nabla \cdot {\bf v}_h \Big) \nn\\
&\qquad\qquad\qquad\qquad 
+\Big(P^n-\Pi^\Gamma_{h}P^n,{\bf v}_h\cdot{\bf n}\Big)_{\partial\Omega}+\sum_{m=1}^M\Big(P^n-\Pi^\Gamma_{h}P^n,[{\bf v}_h\cdot{\bf n}]\Big)_{\Gamma_m},
\label{erre-FEM-1}\\[3pt]
& \Big(\nabla\cdot( \U^{n}_h- Q_h\U^{n}) ,\, \varphi_h\Big) =0,
\label{erre-FEM-2}\\[3pt]
& \biggl(D( \U^{n}_h,x)^{-1}{\bf W}_h^{n+1}-D( \U^{n},x)^{-1}\overline{\bf W}_h^{n+1},\,\overline {\bf v}_h\biggl) =\Big({\cal C}^{n+1}_h-\overline{\cal C}^{n+1}_h,\,
\nabla \cdot \overline {\bf v}_h \Big) \nn\\
&\qquad\qquad\qquad\qquad
+\Big({\cal C}^{n+1}-\Pi^\Gamma_{h}{\cal C}^{n+1},\overline {\bf v}_h\cdot{\bf n}\Big)_{\partial\Omega}
+\sum_{m=1}^M\Big({\cal C}^{n+1}-\Pi^\Gamma_{h}{\cal C}^{n+1},[\overline {\bf v}_h\cdot{\bf n}]\Big)_{\Gamma_m} ,
\label{erre-FEM-3}\\[3pt]
& \Big(\Phi D_\tau({\cal C}^{n+1}_h-{\cal C}^{n+1}), \, \overline\varphi_h\Big) +
\Big(\nabla\cdot({\bf W}_h^{n+1}-\overline{\bf W}_h^{n+1}),
\, \overline\varphi_h \Big) \nn\\
&~~~-\Big(D(\U^n_h,x)^{-1} ({\bf W}_h^{n+1}-{\bf W}^{n+1})\cdot\U_h^{n}, \, \overline \varphi_h\Big) \label{erre-FEM-4}\\
&   - \Big(\big(D(\U^n_h,x)^{-1}  \U^{n}_h-D(\U^n,x)^{-1}  \U^{n}\big)\cdot{\bf W}^{n+1},\, \overline \varphi_h\Big)
 +\Big(({\cal C}^{n+1}_h-{\cal C}^{n+1}) q_I^{n+1}, \,
\overline\varphi_h\Big) =0. \nn
\end{align}

Firstly,
we take
${\bf v}_h=\U^{n}_h-Q_h \U^{n}$ in (\ref{erre-FEM-1}) and get
\begin{align*}
&\biggl(\frac{\mu({\cal C}^{n}_h)}{k(x)}\big( \U^{n}_h-Q_h
\U^{n}\big) +\frac{\mu({\cal C}^{n}_h)}{k(x)}\big(Q_h \U^{n}-
\U^{n}\big) +\frac{\mu({\cal C}^{n}_h)-\mu({\cal C}^{n})}{k(x)}  \U^{n}~ ,  ~\U^{n}_h -Q_h \U^{n}\biggl)\\
&=\Big(P^n-\Pi^\Gamma_{h}P^n, (\U^{n}_h-Q_h
\U^{n})\cdot{\bf n}\Big)_{\partial\Omega}+\sum_{m=1}^M\Big(P^n-\Pi^\Gamma_{h}P^n,[( \U^{n}_h-Q_h
\U^{n})\cdot{\bf n}]\Big)_{\Gamma_m}\\
&\leq C\|P^n\|_{\overline H^3}\| \U^{n}_h-Q_h
\U^{n}\|_{L^2}h^{2}
\end{align*}
which implies that
\begin{align}\label{UitCh}
\big\| \U^{n}_h - \U^{n}\big\|_{L^2} &\leq C(h^2+\big\|{\cal C}^{n}_h-{\cal C}^{n}\big\|_{L^2}),\quad\mbox{for}~~ n=0,1,\cdots,N .
\end{align}
In particular, we have
\begin{align*}
\big\| \U^{0}_h - L_h\U^{0}\big\|_{L^2} &\leq C(h^2+\big\|{\cal C}^{0}_h-{\cal C}^{0}\big\|_{L^2})\leq Ch^2   .
\end{align*}
and so, by the inverse inequality,
\begin{align*}
\big\| \U^{0}_h - L_h\U^{0}\big\|_{L^\infty} &\leq Ch^{-d/2}\big\| \U^{0}_h - L_h\U^{0}\big\|_{L^2}\leq Ch^{2-d/2}.
\end{align*}
As a result, there exists a positive constant $h_2$ such that when $h<h_2$ we have
\begin{align}
\big\| \U^{0}_h\big\|_{L^\infty}  \leq \| \U^{0}\|_{L^\infty}+1.
\end{align}

Secondly, we proceed with a mathematical induction on
\begin{align}\label{UitCh2}
\|\U^n_h\|_{L^\infty}\leq \|\U^n\|_{L^\infty}+1 ,
\end{align}
which is already proved for $n=0$. In the following, we assume that it holds for $0\leq n\leq k$ and try to prove that it also holds for $n=k+1$.

Taking $\overline {\bf v}_h={\bf W}_h^{n+1}-\overline{\bf W}_h^{n+1}$ and $\overline\varphi_h={\cal C}_h^{n+1}-\overline {\cal C}_h^{n+1}$ in \refe{erre-FEM-3}-\refe{erre-FEM-4}, we obtain
\begin{align*}
& D_\tau\bigg(\frac{1}{2}\|\sqrt{\Phi}({\cal C}^{n+1}_h-\overline{\cal C}^{n+1}_h)\|_{L^2}^2\bigg) +\Big(D( \U^{n}_h,x)^{-1}({\bf W}_h^{n+1}-\overline{\bf W}_h^{n+1}),{\bf W}_h^{n+1}-\overline{\bf W}_h^{n+1}\Big)\\
&=-\Big((D( \U^{n}_h,x)^{-1}-D( \U^{n},x)^{-1})\overline{\bf W}_h^{n+1},\,{\bf W}_h^{n+1}-\overline{\bf W}_h^{n+1}\Big) \\
&~~~+\Big(D(\U^n_h,x)^{-1}\U_h^{n} \cdot({\bf W}_h^{n+1}-{\bf W}^{n+1}), \, {\cal C}^{n+1}_h-\overline{\cal C}^{n+1}_h\Big)\\
&~~~+ \Big(\big(D(\U^n_h,x)^{-1}  \U^{n}_h-D(\U^n,x)^{-1}  \U^{n}\big)\cdot{\bf W}^{n+1},\, {\cal C}^{n+1}_h-\overline{\cal C}^{n+1}_h\Big)\\
&~~~-\Big(({\cal C}^{n+1}_h- {\cal C}^{n+1}) q_I^{n+1}, \,
{\cal C}^{n+1}_h-\overline{\cal C}^{n+1}_h\Big) 
+ \Big(\Phi D_\tau({\cal C}^{n+1}-\overline{\cal C}^{n+1}_h), \, {\cal C}_h^{n+1}-\overline {\cal C}_h^{n+1}\Big)\\
&~~~+\Big({\cal C}^{n+1}-\Pi^\Gamma_{h}{\cal C}^{n+1},
\overline {\bf v}_h\cdot{\bf n}\Big)_{\partial\Omega}
+\sum_{m=1}^M\Big({\cal C}^{n+1}-\Pi^\Gamma_{h}{\cal C}^{n+1},
[\overline {\bf v}_h\cdot{\bf n}]\Big)_{\Gamma_m}\\
&\leq C\|\overline{\bf W}_h^{n+1}\|_{L^\infty}
\|\U^n_h-\U^n\|_{L^2}\|{\bf W}_h^{n+1}-\overline{\bf W}_h^{n+1}\|_{L^2} \\
&~~~+C\|D(\U^n_h,x)^{-1}\U_h^{n}\|_{L^\infty}
(\|{\bf W}_h^{n+1}-\overline{\bf W}_h^{n+1}\|_{L^2} 
+\|{\bf W}^{n+1}-\overline{\bf W}_h^{n+1}\|_{L^2}) 
\|{\cal C}^{n+1}_h-\overline{\cal C}^{n+1}_h\|_{L^2}\\
&~~~+C\|{\bf W}^{n+1}\|_{L^\infty}\|\U^n_h-\U^n\|_{L^2}
\|{\cal C}^{n+1}_h-\overline{\cal C}^{n+1}_h\|_{L^2}
+C\|{\cal C}^{n+1}_h-\overline{\cal C}^{n+1}_h\|_{L^2}
\|{\cal C}^{n+1}-\overline{\cal C}^{n+1}_h\|_{L^2}\\
&~~~+C\|D_\tau({\cal C}^{n+1}-\overline{\cal C}^{n+1}_h)\|_{L^2}
\| {\cal C}_h^{n+1}-\overline {\cal C}_h^{n+1}\|_{L^2}
+C\|{\cal C}^{n+1}\|_{\overline H^3}
\|{\bf W}_h^{n+1}-\overline{\bf W}_h^{n+1}\|_{L^2}h^2\\
&\leq \frac{1}{2}\Big(D( \U^{n}_h,x)^{-1}
({\bf W}_h^{n+1}-\overline{\bf W}_h^{n+1}),{\bf W}_h^{n+1}-\overline{\bf W}_h^{n+1}\Big) \\
&~~~+C(1+\|\U^n_h\|_{L^\infty})(
\|{\cal C}^{n+1}_h-\overline{\cal C}^{n+1}_h\|_{L^2}^2
+\|{\bf  U}^{n}_h-{\bf  U}^{n}\|_{L^2}^2
+\|{\bf W}^{n+1}-\overline{\bf W}_h^{n+1}\|_{L^2}^2
+\|{\cal C}^{n+1}-\overline{\cal C}^{n+1}_h\|_{L^2}^2)  \\
&\leq \frac{1}{2}\Big(D( \U^{n}_h,x)^{-1}
({\bf W}_h^{n+1}-\overline{\bf W}_h^{n+1}),
{\bf W}_h^{n+1}-\overline{\bf W}_h^{n+1}\Big) \\
&~~~+C(1+\|\U^n_h\|_{L^\infty})
(\|{\cal C}^{n+1}_h-\overline{\cal C}^{n+1}_h\|_{L^2}^2
+\|{\cal C}^{n}_h-\overline{\cal C}^{n}_h\|_{L^2}^2
+\|D_\tau({\cal C}^{n+1}-\overline{\cal C}^{n+1}_h)\|_{L^2}^2+h^4 )  .
\end{align*}

Since $\displaystyle \|{\bf U}^n\|_{L^\infty}\leq C \|{\bf U}^n\|_{\overline H^2}\leq C$, by applying \refe{UitCh2} and Gronwall's inequality, there exists a positive constant $\tau_4$ such that when $\tau<\tau_4$ we have
\begin{align*}
\max_{0\leq n\leq k}\|{\cal C}^{n+1}_h-\overline{\cal C}^{n+1}_h \|_{L^2}^2
&\leq Ch^2 \, ,
\end{align*}
which, together with \refe{p-app}, gives
\begin{align}
& \max_{0\leq n\leq k}\|{\cal C}^{n+1}_h-{\cal C}^{n+1}  \|_{L^2}\leq Ch^2 .
\end{align}
From \refe{UitCh} we further derive that
\begin{align}
\max_{0\leq n\leq k}\big\| \U^{n+1}_h - \U^{n+1}\big\|_{L^2} \leq Ch^2 ,
\end{align}
which implies that
\begin{align*}
\big\| \U^{k+1}_h - L_h\U^{k+1}\big\|_{L^2} &\leq Ch^2,
\end{align*}
and so, by the inverse inequality,
\begin{align*}
\big\| \U^{k+1}_h - L_h\U^{k+1}\big\|_{L^\infty} &\leq Ch^{-d/2}\big\| \U^{k+1}_h - L_h\U^{k+1}\big\|_{L^2}\leq Ch^{2-d/2}.
\end{align*}
In view of the last inequality, there exists a positive constant $h_3$ such that when $h<h_3$ we have
\begin{align}
\big\| \U^{k+1}_h\big\|_{L^\infty}  \leq \| \U^{k+1}\|_{L^\infty}+1.
\end{align}
The mathematical induction on \refe{UitCh2} is completed, and Proposition \ref{BdUnh} is proved with $\tau_*=\min(\tau_0,\tau_4)$ and $h_*:=\min(h_1,h_2,h_3)$. ~
 \endproof

\section{Proof of Theorem \ref{MainTHM}}
\label{SEction5}
\setcounter{equation}{0}
In this section, we prove Theorem \ref{MainTHM} 
based on the boundedness 
of the fully discrete solution proved in Proposition \ref{BdUnh}.

Similar as the last section, for any fixed integer $n\geq {-1}$
we introduce the Ritz projection
 $(\overline c^{n+1}_h,\overline {\bf w}_h^{n+1})\in
S_h^r\times {\bf H}^r_h$ of $(c^{n+1},{\bf w}^{n+1})\in 
H^1\times {\bf H}_\Gamma^1$ as the finite element 
solution of
\begin{align}\label{MEPrjc3}
\left\{
\begin{array}{ll}
(\nabla\cdot (\overline {\bf w}_h^{n+1}-{\bf w}^{n+1}),\chi_h)=0, 
&\forall~\chi_h\in S_h^r, \\[5pt]
(D(\u^{n+1},x)^{-1}(\overline {\bf w}_h^{n+1}-{\bf w}^{n+1}),{\bf v}_h) 
=(\overline c^{n+1}_h-c^{n+1},\nabla\cdot {\bf v}_h),&\forall~ {\bf v}_h\in {\bf H}_h^r ,
\end{array}
\right.
\end{align}
with $\int_\Omega(\overline c^{n+1}_h-c^{n+1})d x=0$ for 
the uniqueness of solution.  
Then there exists a positive 
constant $h_{**}\leq h_*$ such that when $h<h_{**}$ 
the following inequalities hold:
\begin{align} 
&\|\phi-\Pi_h\phi\|_{L^2}\leq
C\|\phi\|_{\overline H^{r+1}}h^{r+1}, \quad \forall~\phi\in \overline H^{r+1},  
\label{p-app33} \\[5pt]
&\max_{0\leq n\leq N}(\|c^{n}-\overline c_h^{n}\|_{L^2}
+\|{\bf u}^{n}-Q_h {\bf u}^{n}\|_{L^2}
+\|{\bf w}^{n}-\overline {\bf w}_h^{n}\|_{L^2} ) \leq
Ch^{r+1} ,\\
&\max_{0\leq n\leq N}( \|{\bf w}^{n}\|_{L^\infty}
+\|\overline {\bf w}_h^{n}\|_{L^\infty}) \leq
C , \label{p-app3}\\
&\bigg(\sum_{n=0}^{N-1}\tau
\|D_\tau (c^{n+1}-\overline c_h^{n+1})\|_{L^2}^2
\bigg)^{1/2}\leq Ch^{r+1} .\label{p-apap3}
\end{align}

Firstly, by choosing ${\bf v}_h={\bf U}^{n}_h$ in \refe{e-FEM-1} and $\varphi_h=P^{n}_h$ in \refe{e-FEM-2}, we derive that 
\begin{align}
&\|\U^{n}_h\|_{L^2}^2\leq C\|P_h^{n}\|_{L^2} .\label{FUL2} 
\end{align}
In order to make use of \refe{FUL2}, we
define $\widetilde g^n$ as the solution of
\begin{align*}
\left\{\begin{array}{ll}
\Delta \widetilde g^n =P_h^{n}
&\mbox{in}~~\Omega_m,\\
\left[\widetilde g^{n}\right]=0,\quad
[\widetilde g^n\cdot{\bf n}]=0 &
\mbox{on}~~\Gamma_m,\\
\nabla \widetilde g^n\cdot{\bf n}=0 &\mbox{on}~~\partial\Omega ,
\end{array}
\right.
\end{align*}
and substitute ${\bf v}_h=-Q_h\big(\nabla \widetilde g^n\big)$ into \refe{e-FEM-1}. Then we obtain
\begin{align*}
\|P_h^{n}\|_{L^2}^2=\Big(\frac{\mu( {\cal C}_h^{n})}{k(x)} \U_h^{n},\,{\bf v}_h \Big)\leq C\|\U_h^{n}\|_{L^2}\|{\bf v}_h\|_{L^2}\leq C\|\U_h^{n}\|_{L^2}\|P_h^{n}\|_{L^2},
\end{align*}
which together with \refe{FUL2} implies that
\begin{align}\label{lda800}
\|P_h^{n}\|_{L^2}+\|\U_h^{n}\|_{L^2}\leq C.
\end{align}
When $h\geq h_{**}$, by the inverse inequality we have
\begin{align}\label{lda8}
\|\U_h^{n}\|_{L^\infty}\leq Ch^{-d/2}\|\U_h^{n}\|_{L^2}
\leq Ch_{**}^{-d/2}\|\U_h^{n}\|_{L^2}
\leq C.
\end{align}
Then we choose 
$\overline {\bf v}_h={\bf W}_h^{n+1}$ in \refe{e-FEM-3} and $\overline\varphi_h={\cal C}_h^{n+1}$ in \refe{e-FEM-4}. With the boundedness of $\|\U^n_h\|_{L^\infty}$, we derive that
\begin{align*}
&D_\tau\bigg(\frac{1}{2}\|\sqrt{\Phi}{\cal C}_h^{n+1}\|_{L^2}^2\bigg)+
\Big(D(\U_h^n,x)^{-1}{\bf W}^{n+1}_h,{\bf W}^{n+1}_h\Big) \\
&\leq
\frac{1}{2}\Big(D(\U_h^n,x)^{-1}{\bf W}^{n+1}_h,{\bf W}^{n+1}_h\Big)+C(1+\|\U^n_h\|_{L^\infty})\|{\cal C}_h^{n+1}\|_{L^2}^2 +C\|\widehat c^{n+1}q_I^{n+1}\|_{L^2}^2 \\
&\leq
\frac{1}{2}\Big(D(\U_h^n,x)^{-1}{\bf W}^{n+1}_h,{\bf W}^{n+1}_h\Big)+C\|{\cal C}_h^{n+1}\|_{L^2}^2 +C\|\widehat c^{n+1}q_I^{n+1}\|_{L^2}^2 .
\end{align*}
Applying Gronwall's inequality, there exists a positive constant $\tau_5<\tau_*$ such that when 
$\tau<\tau_5$ and $h\geq h_{**}$ we have
\begin{align}\label{lda7}
\|{\cal C}_h^{n+1}\|_{L^2}\leq C.
\end{align}

Secondly, we assume that $\tau<\tau_5$ and $h<h_{**}$ so that
Proposition \ref{BdUnh} holds. Note that
the mixed formulation of \refe{e-fuel-1}-\refe{e-fuel-3} gives
\begin{align*}
& \Big(\frac{\mu( c^{n})}{k(x)} \u^{n},\,{\bf v}_h \Big)
=\Big(p^{n} ,\, \nabla \cdot {\bf v}_h  \Big)-\Big(p^n,{\bf v}_h\cdot{\bf n}\Big)_{\partial\Omega}-\sum_{m=1}^M\Big(p^n,[{\bf v}_h\cdot{\bf n}]\Big)_{\Gamma_m},
,
\\[3pt]
& \Big(\nabla\cdot \u^{n} ,\, \varphi_h\Big) =\Big(q_I^{n}-q_P^{n},\,
\varphi_h\Big),
\\[3pt]
& \Big(D( \u^{n+1},x)^{-1}{\bf w}^{n+1},\,\overline {\bf v}_h\Big)
=\Big(c^{n+1} ,\, \nabla \cdot \overline {\bf v}_h  \Big)-\Big(c^{n+1},\overline {\bf v}_h\cdot{\bf n}\Big)_{\partial\Omega}
-\sum_{m=1}^M\Big(c^{n+1},[\overline {\bf v}_h\cdot{\bf n}]\Big)_{\Gamma_m},,
\\[3pt]
& \Big(\Phi(x) D_\tau c^{n+1}, \,\overline\varphi_h\Big)
+ \Big(\nabla\cdot {\bf w}^{n+1}, \, \overline\varphi_h \Big)  
- \Big(D(\u^{n+1},x)^{-1}\u^{n}\cdot {\bf w}^{n+1},\, \overline \varphi_h\Big) \nn\\
&\qquad\qquad\qquad\qquad\qquad\qquad\qquad\qquad = \Big(\hat
c^{n+1}q_I^{n+1}-c^{n+1} q_I^{n+1}, \, \overline\varphi_h\Big) + 
({\cal E}^{n+1},\overline\varphi_h)
\end{align*}
for any ${\bf v}_h,\overline {\bf v}_h\in {\bf H}_h^r$ and 
$\varphi_h,\overline\varphi_h\in S_h$. The
above equations with the finite element system
(\ref{e-FEM-1})-(\ref{e-FEM-4}) imply that
\begin{align}
& \biggl(\frac{\mu({\cal C}^{n}_h)}{k(x)} \U^{n}_h
-\frac{\mu(c^{n})}{k(x)} \u^{n},\,{\bf v}_h\biggl) =\Big(P^{n}_h-\Pi_hp^{n},\,
\nabla \cdot {\bf v}_h \Big)\nn\\
&\qquad\qquad\qquad\qquad 
+\Big(p^n-\Pi^\Gamma_{h}p^n,
{\bf v}_h\cdot{\bf n}\Big)_{\partial\Omega}
+\sum_{m=1}^M\Big(p^n-\Pi^\Gamma_{h}p^n,[{\bf v}_h\cdot{\bf n}]\Big)_{\Gamma_m},
\label{erre-FEM-13}\\[3pt]
& \Big(\nabla\cdot( \U^{n}_h- Q_h\u^{n}) ,\, \varphi_h\Big) =0,
\label{erre-FEM-23}\\[3pt]
& \biggl(D( \U^{n}_h,x)^{-1}{\bf W}_h^{n+1}
-D( \u^{n+1},x)^{-1}\overline {\bf w}_h^{n+1},
\,\overline {\bf v}_h\biggl) 
=\Big({\cal C}^{n+1}_h-\overline c^{n+1}_h,\,
\nabla \cdot \overline {\bf v}_h \Big)\nn\\
&\qquad\qquad\qquad\qquad
+\Big(c^{n+1}-\Pi^\Gamma_{h}c^{n+1},
\overline {\bf v}_h\cdot{\bf n}\Big)_{\partial\Omega}
+\sum_{m=1}^M\Big(c^{n+1}-\Pi^\Gamma_{h}c^{n+1},[\overline {\bf v}_h\cdot{\bf n}]\Big)_{\Gamma_m} ,
\label{erre-FEM-33}\\[3pt]
& \Big(\Phi D_\tau({\cal C}^{n+1}_h-c^{n+1}), \, \overline\varphi_h\Big) +
\Big(\nabla\cdot({\bf W}_h^{n+1}-\overline {\bf w}_h^{n+1}),
\, \overline\varphi_h \Big) \nn\\
&~~~-\Big(D(\U^n_h,x)^{-1} ({\bf W}_h^{n+1}-{\bf w}^{n+1})
\cdot\U_h^{n}, \, \overline \varphi_h\Big)  \nn\\
& ~~~ -\Big(\big(D(\U^n_h,x)^{-1}  \U^{n}_h
-D(\u^{n+1},x)^{-1}  \u^{n}\big)\cdot{\bf w}^{n+1},\, \overline \varphi_h\Big)
 +\Big(({\cal C}^{n+1}_h-c^{n+1}) q_I^{n+1}, \,
\overline\varphi_h\Big) \nn \\
&=-({\cal E}^{n+1},\overline\varphi_h).
\label{erre-FEM-43}
\end{align}

Taking
${\bf v}_h=\U^{n}_h-Q_h \u^{n}$ in (\ref{erre-FEM-13}), we get
\begin{align*}
&\biggl(\frac{\mu({\cal C}^{n}_h)}{k(x)}\big( \U^{n}_h-Q_h
\u^{n}\big) +\frac{\mu({\cal C}^{n}_h)}{k(x)}\big(Q_h \u^{n}-
\u^{n}\big)
+\frac{\mu({\cal C}^{n}_h)-\mu(c^{n})}{k(x)}  \u^{n} ,
~  \U^{n}_h -Q_h \u^{n}\biggl)\\
&=\Big(p^n-\Pi^\Gamma_{h}p^n,( \U^{n}_h-Q_h
\u^{n})\cdot{\bf n}\Big)_{\partial\Omega}+\sum_{m=1}^M\Big(p^n-\Pi^\Gamma_{h}p^n,[( \U^{n}_h-Q_h
\u^{n})\cdot{\bf n}]\Big)_{\Gamma_m},\\
&\leq C\|p^n\|_{\overline H^{r+1}}\| \U^{n}_h-Q_h
\u^{n}\|_{L^2}h^{r+1} ,
\end{align*}
which implies that
\begin{align}\label{UitCh3}
\big\| \U^{n}_h - \u^{n}\big\|_{L^2} &\leq C(h^{r+1}
+\big\|{\cal C}^{n}_h-c^{n}\big\|_{L^2}),\quad\mbox{for}~~ n=0,1,\cdots,N .
\end{align}

Taking $\overline {\bf v}_h={\bf W}_h^{n+1}
-\overline {\bf w}_h^{n+1}$ and $\overline\varphi_h={\cal C}_h^{n+1}-\overline c_h^{n+1}$
 in \refe{erre-FEM-33}-\refe{erre-FEM-43}, we obtain
\begin{align*}
& D_\tau\bigg(\frac{1}{2}\|\sqrt{\Phi}({\cal C}^{n+1}_h
-\overline c^{n+1}_h)\|_{L^2}^2\bigg) +\Big(D( \U^{n}_h,x)^{-1}({\bf W}_h^{n+1}
-\overline {\bf w}_h^{n+1}),{\bf W}_h^{n+1}-\overline {\bf w}_h^{n+1}\Big)\\
&=\Big((D( \U^{n}_h,x)^{-1}-D( \u^{n+1},x)^{-1})
\overline {\bf w}_h^{n+1},\,{\bf W}_h^{n+1}-\overline {\bf w}_h^{n+1}\Big) \\
&~~~+\Big(D(\U^n_h,x)^{-1}\U_h^{n} 
\cdot({\bf W}_h^{n+1}-{\bf w}^{n+1}), \, {\cal C}^{n+1}_h-\overline c^{n+1}_h\Big)\\
&~~~+ \Big(\big(D(\U^n_h,x)^{-1}  \U^{n}_h
-D(\u^{n+1},x)^{-1}  \u^{n}\big)\cdot{\bf w}^{n+1},\, {\cal C}^{n+1}_h-\overline c^{n+1}_h\Big)\\
&~~~-\Big(({\cal C}^{n+1}_h- c^{n+1}) q_I^{n+1}, \,
{\cal C}^{n+1}_h-\overline c^{n+1}_h\Big)
+ \Big(\Phi D_\tau(c^{n+1}-\overline c^{n+1}_h), 
\, {\cal C}_h^{n+1}-\overline c_h^{n+1}\Big)\\
&~~~+\Big(c^{n+1}-\Pi^\Gamma_{h}c^{n+1},
\overline {\bf v}_h\cdot{\bf n}\Big)_{\partial\Omega}
+\sum_{m=1}^M\Big(c^{n+1}-\Pi^\Gamma_{h}c^{n+1},
[\overline {\bf v}_h\cdot{\bf n}]\Big)_{\Gamma_m} -\Big({\cal E}^{n+1},
{\cal C}_h^{n+1}-\overline c_h^{n+1}\Big) \\
&\leq C\|\overline {\bf w}_h^{n+1}\|_{L^\infty}
\|\U^n_h-\u^{n+1}\|_{L^2}\|{\bf W}_h^{n+1}-\overline {\bf w}_h^{n+1}\|_{L^2} \\
&~~~+C\|D(\U^n_h,x)^{-1}\U_h^{n}\|_{L^\infty}(\|{\bf W}_h^{n+1}
-\overline {\bf w}_h^{n+1}\|_{L^2} +\|{\bf w}^{n+1}
-\overline {\bf w}_h^{n+1}\|_{L^2}) \|c^{n+1}_h-\overline c^{n+1}_h\|_{L^2}\\
&~~~+C\|{\bf w}^{n+1}\|_{L^\infty}
\|\U^n_h-\u^{n+1}\|_{L^2}\|c^{n+1}_h-\overline c^{n+1}_h\|_{L^2}
+C\|{\cal C}^{n+1}_h-\overline c^{n+1}_h\|_{L^2}^2
+C\|{\cal E}^{n+1}\|_{L^2}^2\\
&\leq \frac{1}{2}\Big(D( \U^{n}_h,x)^{-1}({\bf W}_h^{n+1}
-\overline {\bf w}_h^{n+1}),{\bf W}_h^{n+1}
-\overline {\bf w}_h^{n+1}\Big) +C\|{\cal E}^{n+1}\|_{L^2}^2\\
&~~~+C(1+\|\U^n_h\|_{L^\infty}^2)(\|{\cal C}^{n}_h
-\overline c^{n}_h\|_{L^2}^2+\|{\cal C}^{n+1}_h
-\overline c^{n+1}_h\|_{L^2}^2+\|D_\tau(c^{n+1}
-\overline c^{n+1}_h)\|_{L^2}^2+h^{2r+2}+\tau^2)  .
\end{align*}
Applying Gronwall's inequality, there exists 
a positive constant $\tau_{**}<\tau_5$ such that, 
when $\tau<\tau_{**}$ and $h<h_{**}$, 
Proposition \ref{BdUnh} holds and the last inequality reduces to
\begin{align*}
\max_{1\leq n\leq N}\|{\cal C}^{n}_h-\overline c^{n}_h \|_{L^2}^2
+\sum_{n=0}^{N-1}\tau\|{\bf W}^{n+1}_h
-\overline {\bf w}^{n+1}_h \|_{L^2}^2  \leq C(\tau +h^{ r+1})^2 \, ,
\end{align*}
which together with \refe{UitCh3} implies that
\begin{align}\label{FNest0}
\max_{1\leq n\leq N}\big(\big\|{\cal C}^{n}_h- c^{n} \big\|_{L^2}
+\big\| \U^{n}_h - \u^{n}\big\|_{L^2}
\big)
+\bigg(\sum_{n=1}^{N}\tau\|{\bf W}^{n}_h
-{\bf w}^{n} \|_{L^2}^2\bigg)^{\frac{1}{2}} 
\leq C(\tau+h^{r+1}) .
\end{align}

To estimate $\|P_h^{n}-p^{n}\|_{L^2}$, we define $g^{n}$ as the solution of
\begin{align*}
\left\{\begin{array}{ll}
\Delta g^{n} =P_h^{n}-\Pi_hp^{n}
&\mbox{in}~~\Omega_m,\\
\left[g^{n}\right]=0,\quad
[\nabla g^{n}\cdot{\bf n}]=0 &
\mbox{on}~~\Gamma_m,\\
\nabla g^{n}\cdot{\bf n}=0 &\mbox{on}~~\partial\Omega ,
\end{array}
\right.
\end{align*}
and substitute ${\bf v}_h=-Q_h\big(\nabla g^{n}\big)$ into \refe{erre-FEM-13}. Since $\|{\bf v}_h\|_{L^2}\leq C\|P_h^{n}-p^{n}\|_{L^2}$, it follows that
\begin{align*}
\|P_h^{n}-\Pi_hp^{n}\|_{L^2}^2
&\leq
C(\|{\cal C}_h^{n}-c^{n}\|_{L^2}
+\|\U_h^{n}-\u^{n}\|_{L^2}+C\|p^n\|_{\overline H^{r+1}}h^{r+1})\|{\bf v}_h\|_{L^2}\\
&\leq C(\|{\cal C}_h^{n}-c^{n}\|_{L^2}
+\|\U_h^{n}-\u^{n}\|_{L^2}+C\|p^n\|_{\overline H^{r+1}}h^{r+1})\|P_h^{n}-\Pi_hp^{n}\|_{L^2}\\
&\leq C(\tau+h^{r+1})\|P_h^{n}-\Pi_hp^{n}\|_{L^2} ,
\end{align*}
which gives
\begin{align}\label{FNest1}
\max_{1\leq n\leq N}\|P_h^{n}-\Pi_hp^{n}\|_{L^2}
&\leq C(\tau+h^{r+1}) .
\end{align}

Finally, when $\tau< \tau_{**}$ and $h\geq h_{**}$ 
we see that \refe{lda800}-\refe{lda7} give
\begin{align}\label{FNest2}
&\max_{1\leq n\leq N}(\|P_h^{n}-p^{n}\|_{L^2}
+\|\U_h^{n}-\u^{n}\|_{L^2}+\|{\cal C}_h^{n}-c^{n}\|_{L^2})
+\bigg(\sum_{n=1}^{N}\tau\|{\bf W}^{n}_h
- {\bf w}^{n}  \|_{L^2}^2\bigg)^{\frac{1}{2}} \nn\\
&\leq C\leq \frac{C}{h_{**}^{r+1}}(\tau+h^{r+1}).
\end{align}
From \refe{FNest0}-\refe{FNest2} we see that Theorem \ref{MainTHM} holds.~\endproof

\section{Proof of Lemma \ref{LemHk0} and Lemma \ref{LemHkP1}}\label{LemHkP000}
\setcounter{proposition}{0} 
\setcounter{lemma}{0} 
\setcounter{equation}{0}

In this section, we prove Lemma \ref{LemHk0} and Lemma \ref{LemHkP1},
which were used in Section \ref{SEction3} to prove the 
uniform piecewise regularity of the solution of the linearized PDEs. 
We shall use the notation $x=(x',x_d)$, with $x'=(x_1,\cdots,x_{d-1})$.

\subsection{Proof of Lemma \ref{LemHk0}}\label{PrfLHk0}
Before we prove Lemma \ref{LemHk0}, we need to introduce some lemmas below.

\begin{lemma}\label{LemAPP}  
{\it Let $S_{R}=\{x\in\R^d: |x'|<R ~\mbox{and}~|x_d|<R\}$, 
$S_R^+=S_R\cap\{x\in\R^d:0<x_d<R\}$, $S_R^-=S_R\cap\{x\in\R^d:-R<x_d<0\}$. Let 
$\Gamma=\{x\in\R^d: x_d=0\}$ and 
$\Gamma_R:=S_R\cap\Gamma$. Suppose that 
$A_{ij}=A_{ji}\in H^{2}(S_{2R}^+)\cap H^{2}(S_{2R}^-)$ 
satisfies the strong ellipticity condition
$$K^{-1}|\xi|^2\leq \sum_{i,j=1}^dA_{ij}(x)\xi_i\xi_j\leq K|\xi|^2 
\quad\mbox{for~ $x\in S_{2R}\backslash\Gamma$~ 
and~ $\xi\in\R^d$} ,$$ 
and $\phi\in H^{3}(S_{2R}^+)
\cap H^{3}(S_{2R}^-)$ is a solution of
\begin{align}\label{TrEllP}
\left\{
\begin{array}{ll}
\displaystyle
-\nabla\cdot\big(A\nabla \phi\big)=f &\mbox{in}~~S_{2R}\backslash\Gamma_{2R},\\[10pt]
\displaystyle [\phi]=0,\quad \big[A\nabla \phi\cdot{\bf n}\big]=g &\mbox{on}~~\Gamma_{2R} .
\end{array}
\right.
\end{align}
Then  
\begin{align}
&\|\phi\|_{\overline H^2(S_R)}\leq C_R\big(\|f\|_{ L^2(S_{3R/2})} 
+ \|g\|_{H^{1/2}(\Gamma_{3R/2})}+\|\phi\|_{\overline W^{1,3}(S_{3R/2})}\big) , \label{LemAPP1}\\
&\|\phi\|_{\overline H^3(S_R)}\leq C_R\big(\|f\|_{\overline H^{1}(S_{2R})} 
+ \|g\|_{H^{3/2}(\Gamma_{2R})}+\|\phi\|_{\overline W^{2,4}(S_{2R})} \big) .
\label{LemAPP2}  
\end{align} 
where
$\|\psi\|_{\overline H^k(S_R)}:=\|\psi\|_{ H^k(S_R^+)}
+\|\psi\|_{ H^k(S_R^-)}$
for any $\psi\in H^k(S_R^+)\cap H^k(S_R^-)$ and nonnegative integer $k$.
 }
\end{lemma}

\noindent{\it Proof}~~~ 
To simplify the notations, we relax the dependence on $R$ 
in the generic constant, and
set $(f_1,f_2)=\int_\Omega f_1(x)f_2(x)\d x$,
$(g_1,g_2)_\Gamma=\int_\Gamma g_1(x')g_2(x')\d x'$.

Differentiating the equation \refe{TrEllP} with respect to 
$x_j$ for some fixed $1\leq j\leq d-1$ and denote 
$\phi_j=\partial_j\phi$, we obtain that
\begin{align}\label{TrEllP1}
\left\{
\begin{array}{ll}
\displaystyle
-\nabla\cdot\big(A\nabla \phi_j\big)=\partial_jf
+\nabla\cdot\big(\partial_jA\nabla \phi\big) 
&\mbox{in}~~S_{2R}\backslash\Gamma_{2R},\\[10pt]
\displaystyle [\phi_j]=0,\quad 
\big[A\nabla \phi_j\cdot{\bf n}\big]
=\partial_jg-\big[\partial_jA\nabla \phi\cdot{\bf n}\big] 
&\mbox{on}~~\Gamma_{2R} .
\end{array}
\right.
\end{align}
where ${\bf n}$ denote the upward unit normal vector on $\Gamma$.

If we define $\zeta_R$ as a smooth cut-off function satisfying $0\leq \zeta_R\leq 1$, 
$\zeta_R=1$ in $S_R$ and $\zeta_R=0$ outside $S_{3R/2}$, then
\refe{TrEllP1} times $\phi_j\zeta_R^2$ gives
\begin{align*}
&(\zeta_R^2 A\nabla \phi_j,\nabla\phi_j )
+(2\zeta_R\phi_j  A\nabla \phi_j,\nabla\zeta_R)  \\
&=-(f\zeta_R^2,\partial_j\phi_j) -(2\phi_j f\zeta_R,\partial_j\zeta_R)
-(\zeta_R^2\partial_jA\nabla \phi,\nabla\phi_j )
-(2\zeta_R\phi_j\partial_jA\nabla \phi,\nabla\zeta_R)+(\partial_jg  ,\phi_j\zeta_R^2)_{\Gamma}  ,
\end{align*}
which reduces to
\begin{align*} 
\|\nabla \phi_j\|_{L^2(S_R)}^2
&\leq C\big(\|\phi_j\|_{L^2(S_{3R/2})}^2+\|f\|_{L^2(S_{3R/2})}^2
+\|\partial_jA\|_{L^6(S_{2R})}^2\|\nabla\phi\|_{L^3(S_{3R/2})}^2  \\
&~~~~~~~~
+\|\partial_jA\|_{L^6(S_{2R})}
\|\nabla \phi\|_{L^{12/5}(S_{3R/2})}^2
+\|g\partial_j\zeta_R\|_{L^2(\Gamma)}
\|\phi_j\zeta_R\|_{L^2(\Gamma)}\\
&~~~~~~~~ +\|\partial_j(g\zeta_R)\|_{H^{-1/2}(\Gamma)}
\|\phi_j\zeta_R\|_{H^{1/2}(\Gamma)}\big) \\
&\leq C\big(\|f\|_{L^2(S_{3R/2})}^2 
+\|g\|_{H^{1/2}(\Gamma_{3R/2})}^2+
\|\phi_j \|_{H^{1}(S_{3R/2})}^2
+\|\nabla\phi\|_{L^3(S_{3R/2})}^2 \big) ,
\end{align*}
and from \refe{TrEllP} we see that
\begin{align*} 
\|\partial_{dd}\phi\|_{\overline L^2(S_R)}&
= \bigg\|A_{dd}^{-1}\bigg(\sum_{(i,j)\neq (d,d)}
A_{ij}\partial_{ij}\phi +\sum_{i,j=1}^d\partial_{i}
A_{ij}\partial_{j}\phi+f\bigg)\bigg\|_{\overline L^2(S_R)} \\
&\leq C\big(\|f\|_{L^2(S_{3R/2})}  
+\|g\|_{H^{1/2}(\Gamma_{3R/2})}+
\|\phi_j\|_{H^{1}(S_{3R/2})}
+\|\nabla\phi\|_{L^3(S_{3R/2})}  \big) .
\end{align*}
The last two inequalities imply \refe{LemAPP1}.

By applying \refe{LemAPP1} to the problem \refe{TrEllP1}, we derive that
\begin{align*} 
\|\phi_j\|_{\overline H^2(S_R)} 
&\leq C\big(\|\partial_jf\|_{\overline L^2(S_{3R/2})} 
+\|\nabla\cdot(\partial_jA\nabla\phi)\|_{\overline L^2(S_{3R/2})}  
+\|\partial_jg\|_{H^{1/2}(\Gamma_{3R/2})}\\
&~~~~~~~ +\|[\partial_jA\nabla\phi\cdot{\bf n}]\|_{H^{1/2}(\Gamma_{3R/2})}
+\|\phi_j\|_{\overline W^{1,3}(S_{3R/2})} \big)\\
&\leq C\big(\|f\|_{\overline H^1(S_{2R})} 
+\|A\|_{\overline H^2(S_{2R})}\|\phi \|_{\overline W^{1,\infty}(S_{2R})}
+\|A\|_{\overline W^{1,6}(S_{2R})}\|\phi \|_{\overline W^{2,3}(S_{2R})}  \\
&~~~~~~~ +\|g\|_{H^{3/2}(\Gamma_{2R})} 
+\|\partial_jA\nabla\phi\cdot{\bf n}\|_{\overline H^{1} (S_{2R})}
+\|\phi_j\|_{\overline W^{1,3}(S_{2R})} \big)\\
&\leq C\big(\|f\|_{\overline H^1(S_{2R})} +\|g\|_{H^{3/2}(\Gamma_{2R})}
+ \|\phi \|_{\overline W^{2,3}(S_{2R})}   \\
&~~~~~~~ +\|A\|_{\overline H^2(S_{2R})}\|\phi\|_{\overline W^{1,\infty} (S_{2R})}
+\|A\|_{\overline W^{1,6}(S_{2R})}\| \phi\|_{\overline W^{2,3} (S_{2R})}
  \big)\\
 &\leq C\big(\|f\|_{\overline H^1(S_{2R})} +\|g\|_{H^{3/2}(\Gamma_{2R})}
+ \|\phi \|_{\overline W^{2,4}(S_{2R})}  \big) .\end{align*}
Then from \refe{TrEllP} we derive that
\begin{align*} 
\|\partial_{dd}\phi_d\|_{\overline L^2(S_R)}&
= \bigg\|A_{dd}^{-1}\bigg(\sum_{(l,k)\neq (d,d)}A_{lk}\partial_{lk}\phi_d 
+\sum_{l,k=1}^d\partial_{l}A_{lk}\partial_{k}\phi_d+\partial_df
+\nabla\cdot\big(\partial_dA\nabla \phi\big) \bigg)\bigg\|_{\overline L^2(S_R)} \\
&\leq C\big(\|f\|_{\overline H^1(S_{2R})}
 +\|g\|_{H^{3/2}(\Gamma_{2R})} 
+ \|\phi \|_{\overline W^{2,4}(S_{2R})}  \big) .
\end{align*}
The last two inequalities imply \refe{LemAPP2}, and
the proof of Lemma \ref{LemAPP} is completed.~ \endproof\bigskip

The above lemma can be easily extended to 
the case that $\Gamma(\varphi)$ is a  smooth 
surface defined by $x_d=\varphi(x')$ for some 
smooth function $\varphi:\R^{d-1}\rightarrow\R$. 
\begin{lemma}\label{Lemm73}
{\it Let $S_{R}(\varphi)=\{x\in\R^d: |x'|<R~\mbox{and}
~\varphi(x')-R<x_d<\varphi(x')+R\}$, $S_R^+(\varphi)
=\{x\in\R^d: |x'|<R~\mbox{and}~\varphi(x')<x_d<\varphi(x')+R\}$, 
$S_R^-=\{x\in\R^d: |x'|<R~\mbox{and}~\varphi(x')-R<x_d<\varphi(x')\}$, 
and $\Gamma_R(\varphi)=S_R(\varphi)\cap\{x_d=\varphi(x')\}$. 
Suppose that $A_{ij}=A_{ji}\in  H^{3}(S_{2R}^+(\varphi))\cap 
H^3(S_{2R}^-(\varphi))$ satisfies that 
$$K^{-1}|\xi|^2\leq \sum_{i,j=1}^dA_{ij}(x)\xi_i\xi_j\leq K|\xi|^2
\quad\mbox{for ~$x\in S_{2R}(\varphi)\backslash\Gamma_{2R}(\varphi)$~ 
and~ $\xi\in\R^d$},$$ and assume that $\phi\in H^{3}(S_{2R}^+(\varphi))\cap 
H^3(S_{2R}^-(\varphi))$ is a solution of
\begin{align}\label{TrmMP1}
\left\{
\begin{array}{ll}
\displaystyle
-\nabla\cdot\big(A\nabla \phi\big)=f &\mbox{in}~~
S_{2R}(\varphi)\backslash\Gamma_{2R}(\varphi),\\[10pt]
\displaystyle [\phi]=0,\quad \big[A\nabla \phi\cdot{\bf n}\big]=g 
&\mbox{on}~~\Gamma_{2R}(\varphi) .
\end{array}
\right. 
\end{align}
Then
\begin{align}
&\|\phi\|_{\overline H^2(S_R(\varphi))}
\leq C_R\big(\| f\|_{ L^2(S_{2R}(\varphi))} 
+ \| g\|_{H^{1/2}(\Gamma_{2R}(\varphi))}
+\|\phi\|_{\overline W^{1,3}(S_{2R}(\varphi))} \big) , 
  \label{TrmMP12}\\[5pt]
&\|\phi\|_{\overline H^3(S_R(\varphi))}
\leq C_R\big(\|f\|_{ \overline H^1(S_{2R}(\varphi))}  
+ \| g\|_{H^{3/2}(\Gamma_{2R}(\varphi))}
+\|\phi\|_{\overline W^{2,4}(S_{2R}(\varphi))} \big)  ,
\label{TrmMP13}
\end{align} 
where
$\|\psi\|_{\overline H^k(S_R(\varphi))}
=\|\psi\|_{ H^k(S_R^+(\varphi))}
+\|\psi\|_{ H^k(S_R^-(\varphi))}$
for any $\psi\in H^k(S_R^+(\varphi))\cup H^k(S_R^-(\varphi))$ 
and nonnegative integer $k$.
 }
\end{lemma}
\noindent{\it Proof}~~~
Let $x=\Psi(y)$ denote the coordinates transformation 
$x'=y'$ and $x_d=y_d+\varphi(y')$. 
Under this coordinates transformation,  the 
 problem \refe{TrmMP1} is converted to
 \begin{align}\label{TrmMP2}
\left\{
\begin{array}{ll}
\displaystyle
-\nabla_y\cdot\big(\widetilde A(y)\nabla_y\widetilde \phi\big)
=\widetilde f(y) &\mbox{in}~~S_{2R} 
\backslash\Gamma_{2R}   ,\\[10pt]
\displaystyle [\widetilde\phi]=0,\quad \big[\widetilde A(y)
\nabla_y \widetilde\phi\cdot{\bf n}\big]=\widetilde g(y) 
&\mbox{on}~~\Gamma_{2R}  ,
\end{array}
\right. 
\end{align}
 where $\widetilde\phi(y)=\phi(\Psi(y))$, 
 $\widetilde A(y)=A(\Psi(y))$, $\widetilde f(y)=f(\Psi(y))$ 
 and $\widetilde g(y)=g(\Psi(y))\sqrt{1+|\nabla\varphi(y')|^2}$.
By applying Lemma \ref{LemAPP} to the 
problem \refe{TrmMP2}, we get
\begin{align*}
&\|\widetilde\phi\|_{H^2(S_R^+)}
+\|\widetilde\phi\|_{H^2(S_R^-)}\\
&\leq C_R\big(\|\widetilde f\|_{ L^2(S_{2R})} 
+ \|\widetilde g\|_{H^{1/2}(\Gamma_{2R})}
+\|\widetilde\phi\|_{W^{1,3}(S_{2R}^+)}
+\|\widetilde\phi\|_{W^{1,3}(S_{2R}^-)}\big) , 
\\
&\|\widetilde\phi\|_{H^3(S_R^+)}+\|\widetilde\phi\|_{H^3(S_R^-)}\\
&\leq C_R\big(\|\widetilde f\|_{H^{1}(S_{2R}^+)} 
+\|\widetilde f\|_{H^{1}(S_{2R}^-)} 
+ \|\widetilde g\|_{H^{3/2}(\Gamma_{2R})}
+\|\widetilde\phi\|_{W^{2,4}(S_{2R}^+)} 
+\|\widetilde\phi\|_{W^{2,4}(S_{2R}^-)} \big) . 
\end{align*} 
Transforming back to the $x$-coordinates, 
the last two inequalities imply \refe{TrmMP12}-\refe{TrmMP13}.~ 
\endproof\bigskip

\noindent{\it Proof of Lemma {\rm\ref{LemHk0}}}~~~
Without loss of generality, we can assume that the 
functions $A_{ij}$, $f$ and $g$ are sufficiently smooth
so that the problem \refe{HkEstTPEq} has a piecewise $H^3$ 
solution \cite{BP70}.
If we can prove \refe{HkEstTP} with a 
constant $C_R$ which does not
depend on the extra smoothness of $A_{ij}$, $f$ and $g$,
then a compactness argument gives
\refe{HkEstTP} for the nonsmooth $A_{ij}$, $f$ and $g$ 
under the condition of Lemma \ref{LemHk0}.

First, multiplying the equation \refe{HkEstTPEq} by $\phi$, 
it is easy to derive the basic $H^1$ estimate:
\begin{align*} 
\|\phi\|_{H^1(\Omega)}\leq C\bigg(\|f\|_{L^2(\Omega)} 
+\sum_{m=1}^M\|g\|_{H^{-1/2}(\Gamma_m)}\bigg) \, .
\end{align*}

Secondly, by a ``partition of unity'', there exist a finite number of cylinders $S_{2R,j}\subset\Omega$, $j=1,\cdots,J$, such that $\{S_{R,j}\}_{j=1}^J$ covers $\Gamma_m$, $m=1,\cdots,M$. Moreover, each $S_{2R,j}$ only intersects one interface $\Gamma_m$ and in each $S_{2R,j}$, up to a rotation, the interface $\Gamma_m$ can be expressed as $x_d=\varphi_j(x')$ for some smooth function $\varphi_j:\R^{d-1}\rightarrow\R$. Then, by applying Lemma \ref{Lemm73}, we derive that
\begin{align*}
&\|\phi\|_{\overline H^2(S_{R,j})}\leq C_R\big(\| f\|_{ L^2(S_{2R,j} )} 
+ \| g\|_{H^{1/2}(\Gamma_{m}\cap S_{2R,j} )}+\|\phi\|_{\overline W^{1,3}(S_{2R,j})} \big), \quad j=1,\cdots,J ,\\
&\|\phi\|_{\overline H^3(S_{R,j})}\leq C_R\big(\| f\|_{ \overline H^1(S_{2R,j} )} 
+ \| g\|_{H^{3/2}(\Gamma_{m}\cap S_{2R,j} )}+\|\phi\|_{\overline W^{2,4}(S_{2R,j})} \big), \quad j=1,\cdots,J .
\end{align*}
 
Let $D=\Omega\backslash \cup_{j=1}^JS_{R,j}$. It is well-known that, by the interior estimates of elliptic equations, there hold
\begin{align*}
&\|\phi\|_{ H^2(D)}\leq C_R\big(\| f\|_{ L^2(\Omega)} 
+\|\phi\|_{H^{1 }(\Omega)} \big)  ,\\
&\|\phi\|_{ H^3(D)}\leq C_R\big(\| f\|_{ \overline H^1(\Omega)} 
+\|\phi\|_{H^{1 }(\Omega)} \big) .
\end{align*}
The last four inequalities imply that
\begin{align*}
\|\phi\|_{\overline H^2(\Omega)}
&\leq C_R\bigg(\| f\|_{ L^2(\Omega)} 
+\sum_{m=1}^M\| g\|_{H^{1/2}(\Gamma_{m})}+\|\phi\|_{\overline W^{1,3 }(\Omega)} \bigg) \\
&\leq C_R\bigg(\| f\|_{ L^2(\Omega)} 
+\sum_{m=1}^M\| g\|_{H^{1/2}(\Gamma_{m})}+C_\epsilon\|\phi\|_{H^{1}(\Omega)} +\epsilon\|\phi\|_{\overline H^{2 }(\Omega)} \bigg)  ,\\
\|\phi\|_{\overline H^3(\Omega)}&
\leq C_R\big(\| f\|_{\overline H^1(\Omega)} 
+\sum_{m=1}^M\| g\|_{H^{3/2}(\Gamma_{m})}+\|\phi\|_{\overline W^{2,4 }(\Omega)} \big) \\
&\leq C_R\bigg(\| f\|_{\overline H^1(\Omega)} 
+\sum_{m=1}^M\| g\|_{H^{3/2}(\Gamma_{m})}+C_\epsilon\|\phi\|_{H^{1}(\Omega)} +\epsilon\|\phi\|_{\overline H^{3}(\Omega)} \bigg)   ,
\end{align*}
where $\epsilon$ can be arbitrarily small. 

Finally, by choosing $\epsilon$ small enough and using the basic $H^1$ estimate, the last two inequalities imply \refe{HkEstTP}.~ \endproof


\subsection{Proof of Lemma \ref{LemHkP1}}

Integrating \refe{HkEstPTEq} against $\phi^{n+1}$, 
it is easy to derive that
\begin{align*} 
\max_{0\leq n\leq N-1}
\|\phi^{n+1}\|_{L^2}^2
+\sum_{n=0}^{N-1}\tau\|\nabla\phi^{n+1}\|_{L^2}^2
\leq C\sum_{n=0}^{N-1}\tau\big(\|f^{n+1}\|_{L^2}^2+\|{\bf g}^{n+1}\|_{L^2}^2\big).
\end{align*}
Then, by setting ${\bf g}^0={\bf g}^1$, integrating \refe{HkEstPTEq} against 
$-\frac{1}{\Phi}\nabla\cdot\big(A^{n+1}\nabla \phi^{n+1}-{\bf g}^{n+1}\big)$ gives
\begin{align*} 
&\big(D_\tau\nabla\phi^{n+1},A^{n+1}
\nabla \phi^{n+1}-{\bf g}^{n+1}\big)
+\frac{1}{2}\big\|\Phi^{-1/2}\nabla\cdot
\big(A^{n+1}\nabla \phi^{n+1}\big)\big\|_{\overline L^2}^2\\
&\leq C\|f^{n+1}\|_{\overline L^2}^2
+C\|\nabla\cdot {\bf g}^{n+1}\|_{\overline L^2}^2
\end{align*}
which further reduces to
\begin{align*} 
&D_\tau\big[\big(A^{n+1}\nabla\phi^{n+1},\nabla \phi^{n+1}\big)
-2\big(\nabla\phi^{n+1},{\bf g}^{n+1}\big)\big]
+\big\|\Phi^{-1/2}\nabla\cdot\big(A^{n+1}\nabla \phi^{n+1}\big)\big\|_{\overline L^2}^2\\
&\leq C\|f^{n+1}\|_{\overline L^2}^2+C\|\nabla\cdot {\bf g}^{n+1}\|_{\overline L^2}^2
+C\big(D_\tau A^{n+1}\nabla\phi^{n}\cdot\nabla \phi^{n})
-\big(\nabla\phi^{n},D_\tau{\bf g}^{n+1}\big)\\
&\leq C\|f^{n+1}\|_{\overline L^2}^2+C\|{\bf g}^{n+1}\|_{\overline H^1}^2
+C\|D_\tau A^{n+1}\|_{\overline L^2}\|\nabla\phi^{n}\|_{\overline L^4}^2
+\|\nabla\phi^{n}\|_{\overline L^2}\|D_\tau{\bf g}^{n+1}\|_{\overline L^2}d_{n,0}\\
&\leq C\|f^{n+1}\|_{\overline L^2}^2+C\|{\bf g}^{n+1}\|_{\overline H^1}^2
+C_\epsilon\|\nabla\phi^{n}\|_{\overline L^2}^2
+\epsilon(\|\nabla\phi^{n}\|_{\overline H^1}^2
+\|D_\tau{\bf g}^{n+1}\|_{\overline L^2}^2d_{n,0}) .
\end{align*}
From Lemma \ref{LemHk0} we know that
\begin{align*} 
\|\phi\|_{\overline H^2}&\leq C\big\|\Phi^{-1/2}\nabla\cdot
\big(A^{n+1}\nabla \phi^{n+1}\big)\big\|_{\overline L^2}
+C\|[{\bf g}^{n+1}\cdot{\bf n}]\|_{H^{1/2}(\Gamma)}\\
&\leq C\big\|\Phi^{-1/2}\nabla\cdot
\big(A^{n+1}\nabla \phi^{n+1}\big)\big\|_{\overline L^2}
+C\|{\bf g}^{n+1}\|_{\overline H^1}
\end{align*}
The last two inequalities imply that
\begin{align*} 
&D_\tau\big[\big(A^{n+1}\nabla\phi^{n+1},\nabla \phi^{n+1}\big)
-2\big(\nabla\phi^{n+1},{\bf g}^{n+1}\big)\big]
+C^{-1}\|\phi\|_{\overline H^2}^2\\
&\leq C\|f^{n+1}\|_{\overline L^2}^2+C\|{\bf g}^{n+1}\|_{\overline H^1}^2
+C_\epsilon\|\nabla\phi^{n}\|_{\overline L^2}^2
+\epsilon(\|\nabla\phi^{n}\|_{\overline H^1}^2
+\|D_\tau{\bf g}^{n+1}\|_{\overline L^2}^2d_{n,0}).
\end{align*}
Summing up the last inequality for $n=0,1,\cdots,m$, we obtain
\begin{align*} 
\max_{0\leq n\leq m}\|\nabla\phi^{n+1}\|_{L^2}^2
+\sum_{n=0}^m \tau\|\phi^{n+1}\|_{\overline H^2}^2 
&\leq C_\epsilon \sum_{n=0}^m
\tau\big(\|f^{n+1}\|_{\overline L^2}^2+\|{\bf g}^{n+1}\|_{\overline H^1}^2
+\|\nabla\phi^{n}\|_{\overline L^2}^2\big)\\
&~~~ +\epsilon\sum_{n=0}^m
\tau(\|\nabla\phi^{n}\|_{\overline H^1}^2
+\|D_\tau{\bf g}^{n+1}\|_{\overline L^2}^2d_{n,0})
+C\|{\bf g}^{m+1}\|_{L^2}^2  ,
\end{align*}
which further reduces to \refe{HkEstPT2}.
The proof of Lemma \ref{LemHkP1} is completed.~ \endproof

\section{Numerical examples}
\setcounter{equation}{0}
In this section, we present numerical examples to support our
theoretical error analysis.
The computations are performed 
with the software FreeFEM++ \cite{Hecht}.

We solve the problem 
\begin{align}
&\left\{
\begin{array}{ll}
\displaystyle\Phi(x)\frac{\partial c}{\partial t}-\nabla\cdot(D({\bf u},x)\nabla
c)+{\bf u}\cdot\nabla c= f &\mbox{in}~~\Omega_0\cup \Omega_1,\\[10pt]
\left[c\right]=0,\quad \left[D(\u,x)\nabla c\cdot{\bf n}\right]=0 &\mbox{on}~~\Gamma ,\\[10pt]
D(\u,x)\nabla c\cdot{\bf n}=0&\mbox{on}~~\partial\Omega,\\[8pt]
c(x,0)=c_0(x)~~ &\mbox{for}~~x\in\Omega ,
\end{array}
\right.\\[8pt]
&\left\{
\begin{array}{ll}
\displaystyle  \nabla\cdot\u= g  &\mbox{in}~~\Omega_0\cup \Omega_1,\\[8pt]
\displaystyle\u=-\frac{k(x)}{\mu(c)}\nabla p  &\mbox{in}~~\Omega_0\cup \Omega_1,\\[10pt]
\displaystyle \left[ p \right]=0,\quad \left[\u\cdot{\bf n}\right]=0 &\mbox{on}~~\Gamma ,\\[8pt]
\displaystyle \u\cdot{\bf n}=0    &\mbox{on}~~\partial\Omega ,
\end{array}
\right.
\end{align}
in the unit ball $\Omega=\{(x_1,x_2): |x_1|^2+|x_2|^2<1\} $ which is separated by the interface $$\Gamma_1=\{(x_1,x_2): |x_1-0.3|^2+|x_2|^2=0.3^2\}$$ into two subdomains 
\begin{align*}
&\Omega_0=\{(x_1,x_2): |x_1-0.3|^2+|x_2|^2>0.3^2\}\quad\mbox{and}\quad
\Omega_1=\{(x_1,x_2): |x_1-0.3|^2+|x_2|^2<0.3^2\}.
\end{align*}
For simplicity, we choose $\mu(c)=1/(1+e^{5c})$, $d_0=d_r=d_p=\alpha_1=\alpha_2=1.0$ and choose the permeability and porosity
\begin{align*}
\Phi(x)=\left\{
\begin{array}{ll}
0.6 &\mbox{for}~~x\in\Omega_0 ,\\
0.4 &\mbox{for}~~x\in\Omega_1 , 
\end{array}\right.
\qquad
k(x)=\left\{
\begin{array}{ll}
0.012 &\mbox{for}~~x\in\Omega_0 ,\\
0.008 &\mbox{for}~~x\in\Omega_1 , 
\end{array}\right.
\end{align*}
which are smooth in each subdomain but discontinuous across the interface $\Gamma_1$.

Let $\Omega_K=\{(x_1,x_2):\sqrt{(x_1-0.3)^2+x_2^2}<0.6\}$ so that $\Omega_1\subset\Omega_K\subset\Omega$. The functions $f$, $g$ and the initial data $c_0$ are chosen corresponding to the exact solution
\begin{align*}
&p(x,t)=\left\{
\begin{array}{ll}
100((x_1-0.3)^2+x_2^2-0.09)(0.36-(x_1-0.3)^2-x_2^2)^4/\Phi(x), &\mbox{for}~~ x\in \Omega_K ,\\[5pt]
0 &\mbox{for}~~x\in \Omega\backslash\Omega_K,\\
\end{array}\right.
\\[5pt]
&c(x,t)=0.5+50\,  p(x,t) \cos(0.4\, x_1)\sin(0.4 \, x_2)\sin(4t) ,
\end{align*}
which satisfy the jump conditions on the interface $\Gamma_1$ and the boundary conditions on $\partial\Omega$, while $\nabla p$ and $\nabla c$ are discontinuous across the interface $\Gamma_1$. 

\begin{figure}[pth]
\begin{minipage}[b]{0.38\linewidth}
\centering
\includegraphics[width=\textwidth]{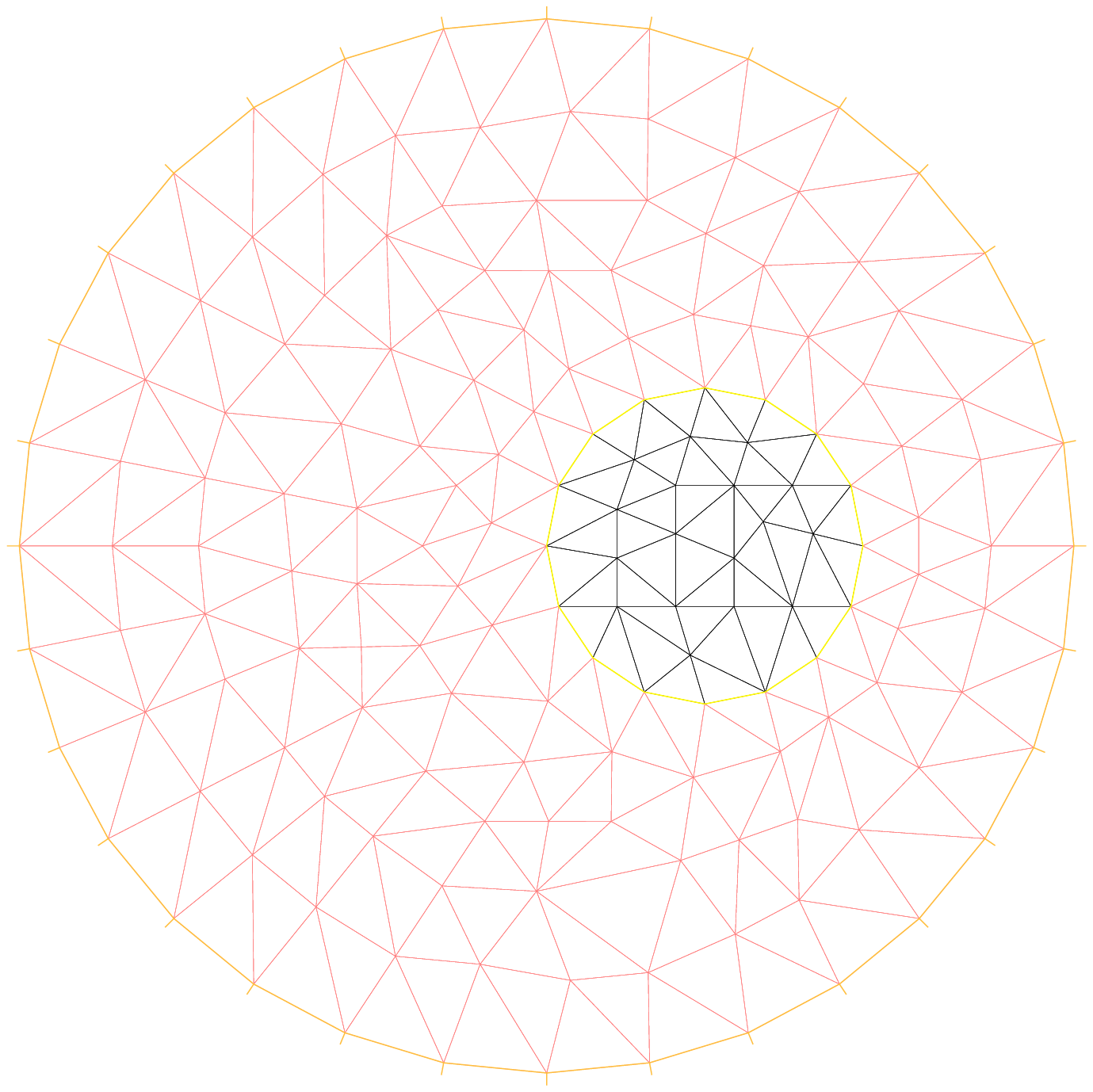}
\end{minipage}
\hspace{-1.8cm}
\begin{minipage}[b]{0.38\linewidth}
\centering
\includegraphics[width=\textwidth]{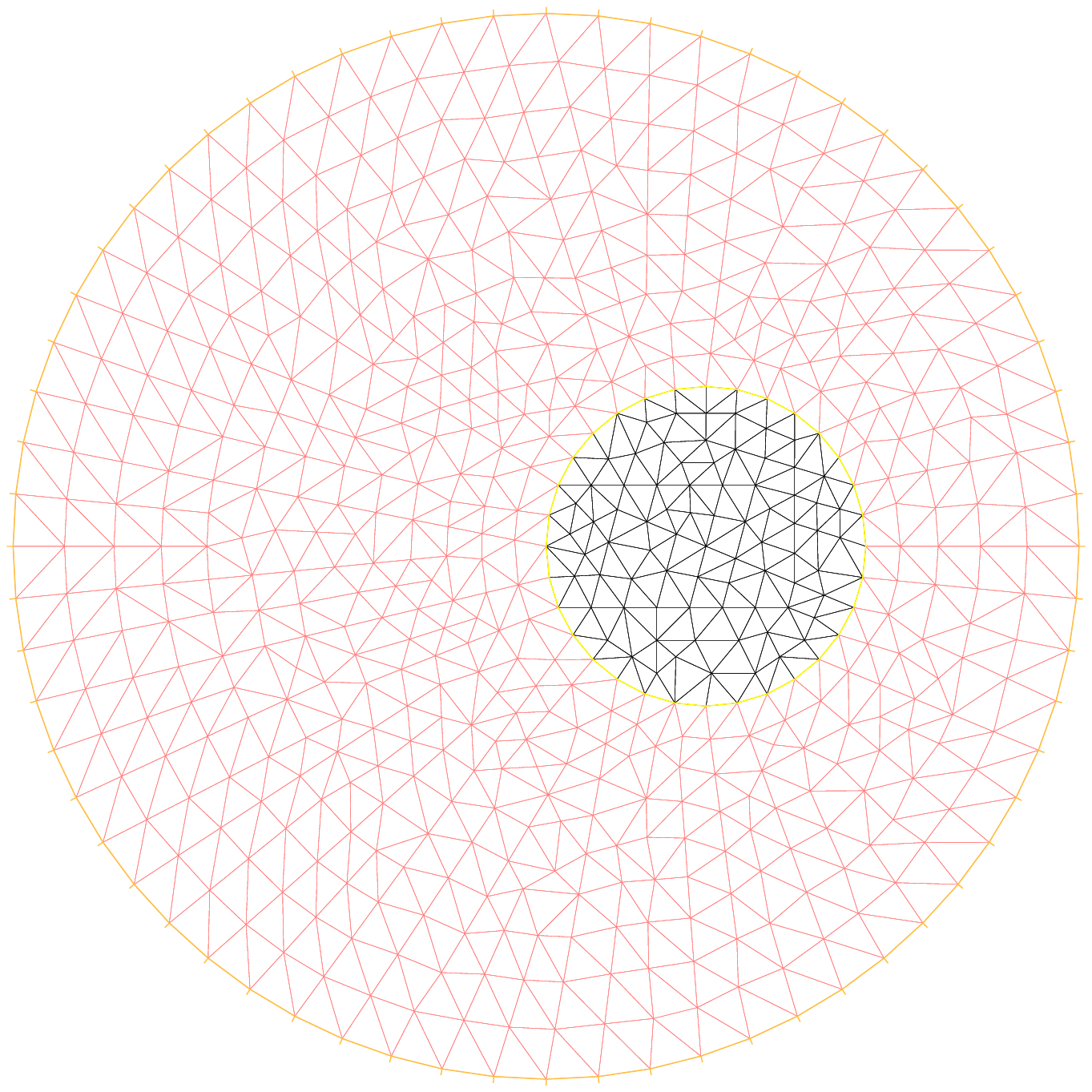}
\end{minipage}
\hspace{-1.8cm}
\begin{minipage}[b]{0.38\linewidth}
\centering
\includegraphics[width=\textwidth]{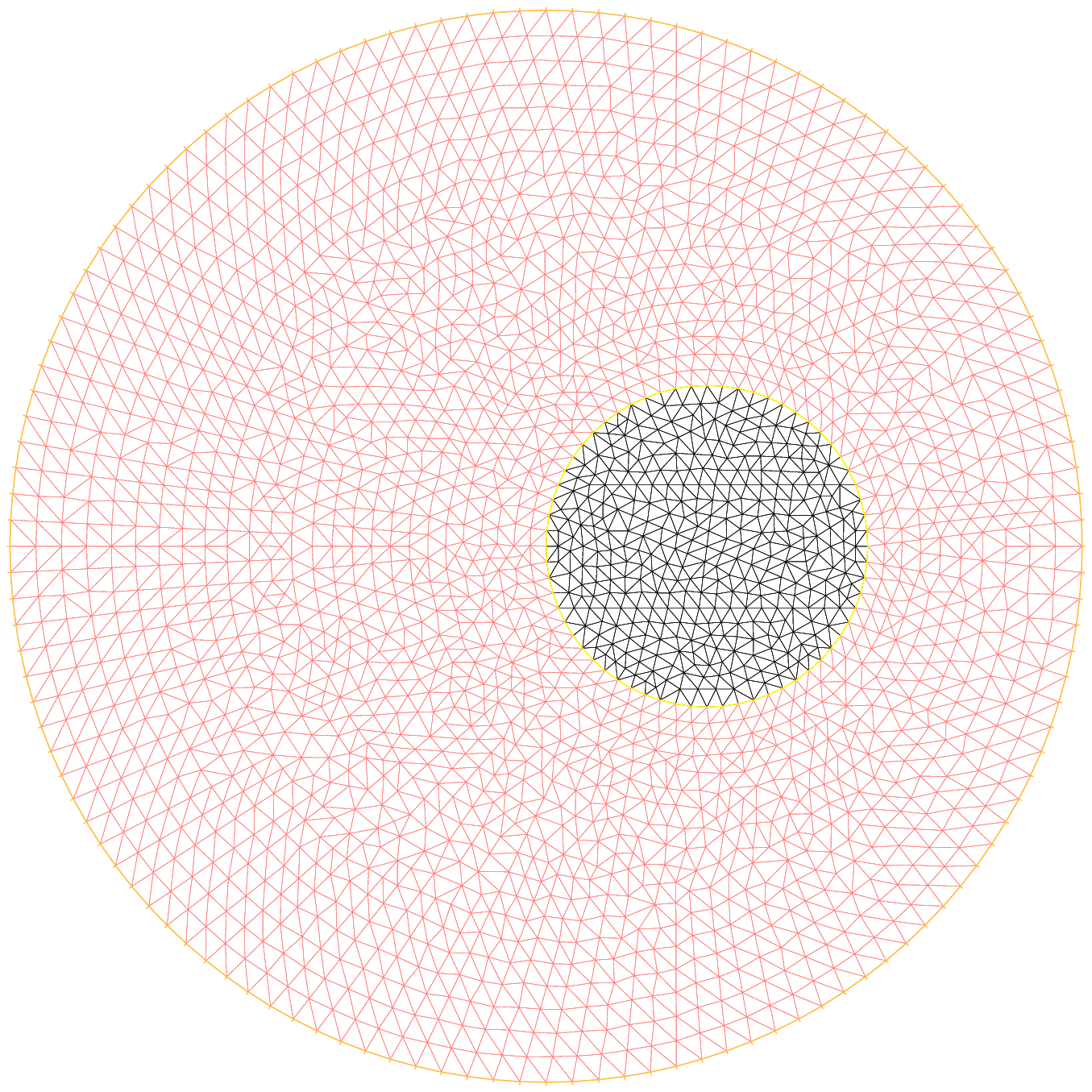}
\end{minipage}
\caption{The finite element meshes 
with $h=1/32$, $h=1/64$ and $h=1/128$.}\label{figure1}
\end{figure}

We partition the domain $\Omega$ into quasi-uniform triangles with $M$ nodes on the boundary $\partial\Omega$ and $M/2$ nodes on the interface $\Gamma$ with $M=32,64,128$, as shown in Figure \ref{figure1}.
For simplicity, we set $h=1/M$ and solve the system by the
proposed method using the Raviart--Thomas linear finite elements up to the time $t=1$.
To test the convergence rate of the proposed method, we solve the problem for different $\tau$ and $h$, 
and present the errors of the numerical solutions 
in Table \ref{ExTab1}, where the convergence rates of ${\bf U}_h^N$ and ${\cal C}^N_h$ are calculated 
by the formulas 
\begin{align*}
&{\rm convergence ~ rate ~ of ~ } {\bf U}_{h}=\ln\big(\|{\bf U}_{h}-{\bf u}\|_{L^2}
/\|{\bf U}_{h/2}-{\bf u}\|_{L^2}\big)/\ln 2 ~,\\
&{\rm convergence ~ rate ~ of ~\, }{\cal C}_{h}\, =\ln\big(\|{\cal C}_{h}-c\|_{L^2}/\|{\cal C}_{h/2}-c\|_{L^2}\big)/\ln 2 ~,
\end{align*}
at the finest two meshes.
From Table \ref{ExTab1} we see that the convergence rate of the numerical solution is about second order, which is consistent with our numerical analysis.

\vskip0.1in

\begin{table}[pth]
\vskip-0.2in
\begin{center}
\caption{Errors of the linearized mixed FEM with $\tau=O(h^2)$.}\vskip 0.1in
\label{ExTab1}
\begin{tabular}{l|l|c|ccc}
\hline
$\tau$  &  $h$
& $\| {\bf U}_h^N - {\bf u}^N \|_{L^2}$
& $\| {\cal C}_h^N - c^N\|_{L^2}$        \\
\hline
1/8   & 1/32
& 3.051E-02  & 1.473E-02 \\
1/32  & 1/64
& 9.769E-03  & 4.280E-03 \\
1/128 & 1/128
& 2.515E-03  & 1.020E-03 \\
\hline
\multicolumn{2}{c|}{convergence rate}
& 1.96       & 2.06 \\
\hline
\end{tabular}
\end{center}
\end{table}

To illustrate the convergence rate with respect to $\tau$, we solve the system for fixed $\tau$ and several different $h$. The errors of the numerical solution are present in Table \ref{ExTab2},
where we can see that the error tends to a constant proportional to $\tau$ 
(as $h$ decreases).

\begin{table}[htp]
\begin{center}
\caption{Errors of the linearized mixed FEM with fixed $\tau$ and refined $h$.}\vskip 0.1in
\label{ExTab2}
\begin{tabular}{l|l|c|ccc}
\hline
$\tau=0.2$ & $h$ 
  & $\| U_h^N - {\bf u}(\cdot ,t_N) \|_{L^2}$ & $\| {\cal C}_h^N - c(\cdot
,t_N) \|_{L^2}$         \\
\cline{2-4}
 & $1/32$  
 & 3.469E-02  & 2.955E-02 \\
 & $1/64$  
 & 2.942E-02  & 3.031E-02 \\
 & $1/96$  
 & 2.887E-02  & 3.025E-02 \\
 & $1/128$  
 & 2.877E-02  & 3.018E-02 \\
\hline
\hline
$\tau=0.1$ & $h$ 
  & $\| U_h^N - {\bf u}(\cdot ,t_N) \|_{L^2}$ & $\| {\cal C}_h^N - c(\cdot
,t_N) \|_{L^2}$         \\
\cline{2-4}
 & $1/32$  
 & 2.840E-02  & 1.438E-02 \\
 & $1/64$  
 & 2.008E-02  & 1.328E-02 \\
 & $1/96$  
 & 1.904E-02  & 1.315E-02 \\
 & $1/128$  
 & 1.883E-02  & 1.308E-02 \\
\hline
\hline
$\tau=0.05$ & $h$ 
  & $\| U_h^N - {\bf u}(\cdot ,t_N) \|_{L^2}$ & $\| {\cal C}_h^N - c(\cdot
,t_N) \|_{L^2}$         \\
\cline{2-4}
 & $1/32$  
 & 2.429E-02  & 1.033E-02 \\
 & $1/64$ 
 & 1.230E-02  & 6.414E-03 \\
 & $1/96$ 
 & 1.029E-02  & 6.100E-03 \\
 & $1/128$ 
 & 9.848E-03  & 6.010E-03 \\
\hline
\end{tabular}
\end{center}
\end{table}

%
%
%



\section{Conclusions}
We have studied the convergence of a 
linearized mixed FEM for a nonlinear elliptic-parabolic interface problem 
from the model of incompressible miscible flow in porous media. 
We showed that the solution of 
the linearized PDEs is piecewise uniformly regular in each 
subdomain separated by the interfaces 
if the solution of the original problem is piecewise
regular, and established optimal-order  
error estimates for the fully discrete solution
without restriction on the grid ratio. 
The analysis presented in this paper, together with 
Lemma \ref{LemHk0}--\ref{LemHkP1}, 
may be extended to other nonlinear parabolic 
interface problems with other time-stepping schemes.

\end{document}